# JUMPING COEFFICIENTS OF MULTIPLIER IDEALS

LAWRENCE EIN, ROBERT LAZARSFELD, KAREN E. SMITH, AND DROR VAROLIN

*Dedicated to Y.-T. Siu on the occasion of his sixtieth birthday*

## Introduction

The purpose of this paper is to study some local invariants attached via multiplier ideals to an effective divisor or ideal sheaf on a smooth complex variety. First considered (at least implicitly) by Libgober [21] and by Loeser and Vaquié [25], [33], [34], these *jumping coefficients* consist of an increasing sequence of positive rational numbers beginning with the log-canonical threshold of the divisor or ideal in question. They encode interesting geometric and algebraic information, and we will see that they arise naturally in several different contexts. Given a polynomial $f \in \mathbf{C}[t_1, \ldots, t_d]$ having only isolated singularities, results of Varchenko, Loeser and Vaquié [36], [25], [33], [34] imply that if $\xi$ is a jumping number of $\{f = 0\}$ lying in the interval $(0, 1]$, then $-\xi$ is a root of the Bernstein-Sato polynomial $b_f(s)$ of $f$. We adapt an argument of Kollár [18] to prove that this holds also when the singular locus of $f$ has positive dimension. In a more algebraic direction, we show that the number of such jumping coefficients bounds the uniform Artin-Rees number of the principal ideal $(f)$ in the sense of Huneke [15]: when $f$ has only isolated singularities, this in turn leads to bounds involving the Milnor number and Tyurina number of $f$. Along the way, we establish a general result relating multiplier to Jacobian ideals. Finally, we discuss the extension of these ideas to the setting of graded families of ideals.

Turning to a more detailed overview, let $X$ be a smooth complex algebraic variety and $x \in X$ a fixed point. Given an effective divisor $A$ on $X$, recall that the *log-canonical threshold* of $A$ at $x$ is the positive (rational) number

$$\mathrm{lct}(A;x) \;=\; \sup \big\{ c \in \mathbf{Q}^{>0} \mid \text{the pair } (X\,,\, c \cdot A) \text{ is log-canonical at } x \big\}.$$

This threshold measures the singularities of $A$ at $x$, with "nastier" singularities corresponding to smaller values of $\mathrm{lct}(A;x)$. It has been the focus of considerable attention in recent years (cf. [18], [6], [27]), and has emerged as a very basic invariant.

Seen from the viewpoint of multiplier ideals, the quantities with which we are concerned are natural generalizations of $\mathrm{lct}(A;x)$. Recall that for every rational $c > 0$ one can define the multiplier ideal sheaf $\mathcal{J}(c \cdot A) = \mathcal{J}(X\,,\, c \cdot A) \subseteq \mathcal{O}_X$. These are

Research of the first author partially supported by NSF Grant DMS 02-00278.
Research of the second author partially supported by NSF Grant DMS 01-39713.
Research of the third author partially supported by NSF Grant DMS 00-70722.

coherent sheaves of ideals, and it is well-known that
$$\mathrm{lct}(A;x) \;=\; \sup\left\{c \in \mathbf{Q}^{>0} \;|\; \mathcal{J}(X,\, c \cdot A)_x = \mathcal{O}_{X,x}.\right\}$$
As $c$ grows larger, the ideals $\mathcal{J}(X,\, c \cdot A)$ decrease, and the idea is simply to look at those (necessarily rational) values of $c$ at which their stalks at $x$ jump. Specifically, it is elementary that there is an increasing discrete sequence
$$0 = \xi_0(A;x) \;<\; \xi_1(A;x) \;<\; \xi_2(A;x) < \ldots$$
of rational numbers $\xi_i = \xi_i(A;x)$ characterized by the properties that
$$\mathcal{J}(X,\, c \cdot A)_x \;=\; \mathcal{J}(X,\, \xi_i \cdot A)_x \quad \text{for} \quad c \in [\xi_i,\, \xi_{i+1}),$$
while
$$\mathcal{J}(X,\, \xi_{i+1} \cdot A)_x \;\subsetneq\; \mathcal{J}(X,\, \xi_i \cdot A)_x$$
for every $i$. (Here we agree by convention that $\mathcal{J}(X,\, 0 \cdot A) = \mathcal{O}_X$.) Thus $\xi_1(A;x) = \mathrm{lct}(A;x)$ is the log-canonical threshold of $A$ at $x$.

**Definition.** The rational numbers $\xi_i(A;x)$ are the *jumping coefficients* or *jumping numbers* of $A$ at $x$. We say that $\xi$ is a jumping coefficient of $A$ on $X$ if it is a jumping number of $A$ at some point $x \in X$.

Given an ideal sheaf $\mathfrak{a} \subseteq \mathcal{O}_X$, the jumping numbers $\xi_i(\mathfrak{a};x) \in \mathbf{Q}$ are defined similarly using the multiplier ideals $\mathcal{J}(X,\, c \cdot \mathfrak{a}) = \mathcal{J}(X,\, \mathfrak{a}^c)$.

We begin in Section 1 by establishing some formal properties of these numbers:

**Proposition A.** *Let $A$ be an effective integral divisor on $X$ passing through $x$.*

  (i). *The collection of all jumping numbers of $A$ is periodic with period $1$. Specifically, $\xi$ is a jumping number for $A$ at $x$ if and only if $(1 + \xi)$ is.*
  (ii). *The jumping numbers of $A$ satisfy the inequality*
(1)
$$\xi_{i+1}(A;x) \;\leq\; \xi_i(A;x) \,+\, \xi_1(A;x)$$
  *for every $i \geq 1$.*

The jumping coefficients associated to an ideal exhibit the analogous periodicity starting at $\dim X - 1$. It follows from (ii) that if the log-canonical threshold $\mathrm{lct}(A;x) = \xi_1(A;x)$ of $A$ at $x$ is small, then $A$ must have many jumping numbers in the interval $[0,1]$. For divisors with isolated singularities, we show that the jumping numbers obey a semi-continuity property analogous to that satisfied by the spectrum of a singularity [32], and we establish a statement of Thom–Sebastiani type for divisors defined by functions in independent sets of variables. Finally, if $X = \mathbf{C}^d$ and $A$ is defined by a polynomial $f \in \mathbf{C}[t_1, \ldots, t_d] = \mathbf{C}[t]$ which is non-degenerate with respect to its Newton polyhedron, then results of Howald [12], [13] permit one to compute the rational numbers $\xi_i(A;0)$ explicitly.

In the case of a polynomial $f \in \mathbf{C}[t] = \mathbf{C}[t_1, \ldots, t_d]$, the jumping numbers of the divisor defined by $f$ are related to some other invariants. To begin with, recall that the *Bernstein-Sato polynomial* $b(s) = b_f(s)$ of $f$ is characterized as the monic polynomial of minimal degree having the property that there is a linear differential operator $P = P(s,t)$ such that
$$Pf^{s+1} = b(s)f^s$$
for all $s$. The Bernstein-Sato polynomial is a very interesting and delicate invariant of the singularities of $\{f = 0\}$ (cf. [26], [17], [24], [18]). A theorem of Kashiwara states that the roots of $b_f(s)$ are negative rational numbers, and Yano [37], Lichten [22] and Kollár [18] show that if $\operatorname{lct}(f)$ is the log-canonical threshold of $f$ on $\mathbf{C}^d$, then $-\operatorname{lct}(f)$ is the largest root $b_f(s)$.[1] We prove that this statement extends in a natural way to higher jumping numbers. Specifically, we use Kollár's argument in [18] to establish:

**Theorem B.** *Let $\xi$ be a jumping coefficient of $f$ on $\mathbf{C}^d$ which lies in the interval $(0,1]$. Then $b_f(-\xi) = 0$.*

When $X$ has an isolated singularity, this is a special case of more precise results due to Varchenko, Loeser and Vaquié ([36], [25], [33], [34]) relating jumping numbers to the spectrum of $f$. In general, not all roots of $b_f(s)$ occur as jumping numbers, but the inequality (1) nonetheless implies the amusing fact that if
$$0 > r_1 > r_2 > \ldots > r_t = -1$$
are the distinct roots of $b_f(s)$ lying in between $0$ and $-1$, then
$$r_{i+1} \geq r_i + r_1$$
for every $1 \leq i < t$. In particular, $r_{i+1} \geq i \cdot r_1 = -i \cdot \operatorname{lct}(A;f)$ for roots in this range. At least in the case of isolated singularities, van Doorn and Steenbrink [35] have established some quite precise inequalities of a closely related nature.

In another direction, we connect jumping coefficients to uniform Artin-Rees numbers. Let $X$ be a smooth complex affine variety of dimension $d$, and $f \in \mathbf{C}[X]$ a non-zero function. Fixing any non-trivial ideal $\mathfrak{b} \subseteq \mathbf{C}[X]$, recall that the Artin-Rees lemma for the principal ideal $(f) \subseteq \mathbf{C}[X]$ states that there is an integer $k > 0$ such that
$$\mathfrak{b}^m \cap (f) \subseteq \mathfrak{b}^{m-k} \cdot (f)$$
for all $m \geq k$. Classically $k$ is allowed to vary with both $f$ and $\mathfrak{b}$, but Huneke [15] proved that one can choose $k$ independent of $\mathfrak{b}$.[2] In this case $k$ is called a uniform Artin-Rees number for $f$, and it is natural to ask what geometric information it depends on.

To formulate our result, we define the *jumping length* $\ell = \ell(f;X)$ of $f$ on $X$ to be the natural number characterized by the property that $\xi_\ell(f;X) = 1$ (it being elementary that $1$ is a jumping number of every principal ideal). In other words, $\ell(f;X)$ counts the number of non-zero jumping coefficients of $f$ which are $\leq 1$. One then has:

---

[1] By definition, the log-canonical threshold $\operatorname{lct}(f)$ of $f$ on $\mathbf{C}^d$ is the minimal value of the log-canonical threshold of the hypersurface defined by $f$ at any of its points. Similarly by the jumping numbers $\xi_i(f)$ of $f$ we mean the jumping numbers of $\operatorname{div}(f)$ on $\mathbf{C}^d$.

[2] The Artin-Rees Lemma and Huneke's theorem actually hold much more generally. Huneke shows in particular that if $\mathfrak{q} \subseteq A$ is an ideal in virtually any Noetherian ring, then there is an integer $k = k(\mathfrak{q})$ depending only on $\mathfrak{q}$ such that $\mathfrak{b}^m \cap \mathfrak{q} \subseteq \mathfrak{q} \cdot \mathfrak{b}^{m-k}$ for every ideal $\mathfrak{b} \subseteq A$ and every $m \geq k$.

**Theorem C.** *In the notation just introduced, $k = d\ell$ is a uniform Artin-Rees number for $f$ at $x$.*

When the hypersurface defined by $f$ has an isolated singularity at some point $x$ but is otherwise smooth, we establish a connection between the jumping length and the Tyurina number $\tau(f;x) = \dim_{\mathbf{C}}\left(\mathcal{O}_x X/(f, \frac{\partial f}{\partial x_1}, \ldots \frac{\partial f}{\partial x_d})\right)$ of $f$ at $x$. Combined with some observations of Huneke in [15], this leads to:

**Corollary D.** *Assume that $(f = 0)$ has an isolated singularity at $x$, and otherwise is smooth. Then $k = \tau(f;x) + d$ is a uniform Artin-Rees number for $f$.*

Huneke had conjectured that there should be a bound of this sort. Along the way to proving Corollary D, we study the connection between multiplier and Jacobian ideals. In this direction we establish a general result of independent interest:

**Theorem E.** *Let $\mathfrak{a} \subseteq \mathcal{O}_X$ be an ideal sheaf on $X$, fix a natural number $m \geq 1$, and denote by $\mathrm{Jac}_m(\mathfrak{a}) \subseteq \mathcal{O}_X$ the $m^{th}$ Jacobian ideal of $\mathfrak{a}$ (Definition 4.1). If the multiplier ideal $\mathcal{J}(\mathfrak{a}^m)$ cuts out an algebraic subset of codimension $m$ in $X$, then*

$$\mathrm{Jac}_m(\mathfrak{a}) \subseteq \mathcal{J}(\mathfrak{a}^{(1-\varepsilon)m}) \quad \text{for all } 0 < \varepsilon \leq 1.$$

The inequality relating jumping lengths to Tyurina numbers comes from taking $\mathfrak{a} = (f)$ and $m = 1$.

Finally, let $\mathfrak{a}_\bullet = \{\mathfrak{a}_k\}$ be a graded system of ideals in the sense of [7]. One can attach jumping coefficients to $\mathfrak{a}_\bullet$ using the asymptotic multiplier ideals $\mathcal{J}(c \cdot \mathfrak{a}_\bullet)$, but now these invariants need no longer be rational or periodic. We give an example to show that the collection $\mathrm{Jump}(\mathfrak{a}_\bullet)$ of jumping numbers can contain cluster points, but we prove that it satisfies the descending chain condition. When all the $\mathfrak{a}_k$ vanish only at a single point the jumping coefficients are discrete, and we show that semicontinuity, as well as the analogue of the inequality (ii) from Proposition A, remain valid in this setting. We also discuss briefly the jumping coefficients attached to a plurisubharmonic function on a complex manifold.

In the case of hypersurfaces, jumping numbers are related to several other constructions that appear in the literature, and we hope that the present paper might call the attention of singularity theorists to multiplier ideals and the invariants they define. Being ourselves at best amateurs in singularity theory, we apologize in advance for any connections or attributions we have overlooked.

The paper is organized as follows. In §1 we define the jumping coefficients, give some examples and establish some basic properties. The connection with Bernstein polynomials and the work of Varchenko appears in §2. Section 3 is devoted to Artin-Rees, where we prove a generalization of Theorem C with $(f)$ replaced by any multiplier ideal. The connection between multiplier and Jacobian ideals appears in §4. Finally, in §5 we give a few results and examples concerning the jumping coefficients attached to graded families of ideals, and discuss briefly the extension to multiplier ideals associated to plurisubharmonic functions. In the hope of convincing the reader that jumping numbers are basic and interesting invariants, we have included numerous concrete examples throughout the paper.

We are grateful to Nero Budur, Igor Dolgachev, Juha Heinonen, Jason Howald, Craig Huneke, János Kollár, Mircea Mustaţă, Toby Stafford, J. Steenbrink and Alex Wolfe for valuable discussions and correspondence. It is a pleasure to dedicate this paper to Y.-T. Siu on the occasion of his sixtieth birthday. Siu's work has played a critical role in developing and promoting the theory of multiplier ideals, and his vision and generosity with ideas have been an inspiration to the authors.

## 1. Definition and Formal Properties

This section is devoted to the definition and basic properties of jumping coefficients.

Let $X$ be a non-singular quasi-projective variety of dimension $d$, and let $\mathfrak{a} \subseteq \mathcal{O}_X$ be an ideal sheaf on $X$. Given any rational (or real) number $c > 0$, recall that one can define the multiplier ideal sheaf
$$\mathcal{J}(X, c \cdot \mathfrak{a}) = \mathcal{J}(X, \mathfrak{a}^c) \subseteq \mathcal{O}_X$$
of $\mathfrak{a}$ with coefficient $c$. This is a coherent sheaf of ideals which roughly speaking measures the singularities of functions $f \in \mathfrak{a}$: for fixed coefficient $c$, "nastier" singularities give rise to "deeper" ideals. These multiplier ideals exhibit remarkable cohomological properties coming from the Kawamata-Viehweg-Nadel vanishing theorem. We refer to [20, Part III] for a detailed development of the theory from an algebro-geometric perspective, to [5, §1] for a quick overview, or to [1] for a gentle introduction to the local picture. The analytic theory appears in [3] and [4].

In practice we assume some familiarity with multiplier ideals, but to fix notation we quickly review the construction. Start by taking any log-resolution $\mu : X' \longrightarrow X$ of $\mathfrak{a}$, so that
$$(2) \qquad \mathfrak{a} \cdot \mathcal{O}_{X'} = \mathcal{O}_{X'}(-F)$$
where $F$ is an effective divisor on $X'$ such that $F + \text{exceptional}(\mu)$ has strict normal crossing support. Writing $K_{X'/X} = K_{X'} - \mu^* K_X$ for the relative canonical bundle of $\mu$, the multiplier ideal in question is defined by setting
$$\mathcal{J}(X, \mathfrak{a}^c) = \mu_* \mathcal{O}_{X'}(K_{X'/X} - [cF]),$$
where $[cF]$ is the integer part (or "round-down") of the $\mathbf{Q}$-divisor (or $\mathbf{R}$-divisor) $cF$. An important point is that the multiplier ideal thus defined is independent of the choice of $\mu$. We sometimes use the notations $\mathcal{J}(c \cdot \mathfrak{a})$ or $\mathcal{J}(\mathfrak{a}^c)$ if the ambient variety is clear, and when $c = 1$ we write simply $\mathcal{J}(\mathfrak{a})$.

**Remark 1.1.** One can similarly attach multiplier ideals to any effective $\mathbf{Q}$-divisor. The construction can be summarized compactly by the equality
$$\mathcal{J}(X, c \cdot A) = \mathcal{J}(X, c \cdot \mathfrak{a})$$

for any integral effective divisor $A$ on $X$ with ideal sheaf $\mathfrak{a} = \mathcal{O}_X(-A)$ and any $c > 0$. The definition for an arbitrary **Q**-divisor $D$ is then determined by the equality $\mathcal{J}(X, c \cdot D) = \mathcal{J}(X, \frac{c}{n} \cdot nD)$ for $n$ a sufficiently divisible integer. □

**Remark 1.2** (Analytic construction of multiplier ideals). Multiplier ideals can be – and originally were – constructed analytically. In fact, choose local generators $f_1 \ldots, f_p \in \mathcal{O}_X$ of $\mathfrak{a}$. Then (the analytic sheaf determined by) $\mathcal{J}(\mathfrak{a}^c)$ is locally generated in a neighborhood of $x$ by all holomorphic functions $g$ such that

$$\frac{|g|^2}{\left(\sum |f_i|^2\right)^c}$$

is integrable near $x$.[3] (See [3] or [4].) □

The basic theme of the present paper is to study the variation of the multiplier ideals $\mathcal{J}(X, \mathfrak{a}^c)$ as a function of $c$. The first point is to analyze the intervals on which these ideals are constant:

**Lemma 1.3.** *With notation as above, let $x \in X$ be a fixed point contained in the zeroes of $\mathfrak{a}$. Then there is an increasing discrete sequence*

$$0 = \xi_0(\mathfrak{a}; x) < \xi_1(\mathfrak{a}; x) < \xi_2(\mathfrak{a}; x) < \ldots$$

*of rational numbers $\xi_i = \xi_i(\mathfrak{a}; x)$ characterized by the properties that*

$$\mathcal{J}(X, c \cdot \mathfrak{a})_x = \mathcal{J}(X, \xi_i \cdot \mathfrak{a})_x \quad \text{for} \quad c \in [\xi_i, \xi_{i+1}),$$

*while*

$$\mathcal{J}(X, \xi_{i+1} \cdot \mathfrak{a})_x \subsetneq \mathcal{J}(X, \xi_i \cdot \mathfrak{a})_x$$

*for every $i$.*

(Here we agree by convention that $\mathcal{J}(X, 0 \cdot \mathfrak{a}) = \mathcal{O}_X$.)

**Definition 1.4** (Jumping coefficients). The rational numbers $\xi_i(\mathfrak{a}; x)$ are the *jumping numbers* or *jumping coefficients* of $\mathfrak{a}$ at $x$. We say that $\xi$ is a jumping coefficient of $\mathfrak{a}$ on a closed subset $Z \subseteq X$ if it is a jumping coefficient of $\mathfrak{a}$ at some point $x \in Z$.

*Proof of Lemma 1.3.* Fix a log resoluton $\mu : X' \longrightarrow X$ of $\mathfrak{a}$ with $\mathfrak{a} \cdot \mathcal{O}_{X'} = \mathcal{O}_{X'}(-F)$, and write

$$F = \sum r_i E_i \quad , \quad K_{X'/X} = \sum b_i E_i,$$

so that

$$\mathcal{J}(X, c \cdot \mathfrak{a}) = \mu_* \mathcal{O}_{X'}\left(\sum (b_j - [cr_j]) \cdot E_j\right).$$

Starting with a given positive rational number $c$, each of the coefficients appearing on the right remains constant if we increase $c$ slightly. Therefore the corresponding multiplier ideals are constant on intervals of the stated shape. Moreover the endpoints of these intervals necessarily occur among those values of $c$ for which

$$ord_{E_j}(K_{X'/X} - cF) = -m$$

---

[3] In other words, $\mathcal{J}(\mathfrak{a}^c)$ is a sheaf of "multipliers" for the weight $\frac{1}{(\sum |f_i|^2)^c}$, hence the name.

for some integer $m \geq 1$ and index $j$ such that $\mu(E_j) \ni x$. In other words the jumping coefficients appear among the numbers $\left\{\frac{b_j+m}{r_j}\right\}$, and therefore the $\xi_i$ are indeed rational and do not have any cluster points. □

**Remark 1.5** (Jumping numbers of **Q**-divisors). The statement of Lemma 1.3 remains valid if the ideal $\mathfrak{a}$ is replaced by an effective **Q**-divisor $D$, and then the jumping numbers of $D$ are defined just as in 1.4. □

**Remark 1.6** (History). The earliest references we know where these invariants appear at least implicitly are the papers [21], [25] of Libgober and Loeser-Vaquié (see especially §IV.3 of [25]); jumping numbers appear more explicitly in the later articles [33, especially p. 1191] and [34, especially pp. 389 - 390] of Vaquié. In this work — which was influenced by the papers [8], [9] of Esnault and Esnault-Viehweg — the invariants in question arose in studying the irregularity of cyclic coverings of $\mathbf{P}^2$ and other surfaces. Most significantly, jumping numbers were related in the cited papers to Hodge-theoretic invariants introduced by Varchenko [36]: this connection is summarized in Remark 2.2. Given that [21], [25], [33] and [34] predate a systematic theory of multiplier ideals, it is perhaps not surprising that what we are calling jumping numbers figure in the cited references largely in passing. □

**Example 1.7** (Integral divisors). If $A$ is an integral divisor then $\mathcal{J}(A) = \mathcal{O}_X(-A)$. But if $c < 1$ then $\mathcal{O}_X(-A)_x \subsetneq \mathcal{J}(c \cdot A)_x$ for every $x \in \mathrm{supp}(A)$ owing to the fact that each of the irreducible components of $A$ evidently appears with smaller multiplicity in the multiplier ideal on the right than in $A$ itself. Therefore $\xi = 1$ is a jumping number of an integral divisor at every point of its support. □

**Example 1.8** (Monomial ideals). A result of Howald [12] allows one to compute explicitly the jumping numbers on $X = \mathbf{C}^d$ of a monomial ideal $\mathfrak{a} \subseteq \mathbf{C}[t_1, \ldots, t_d] = \mathbf{C}[X]$. Identifying monomials with points in $\mathbf{N}^d \subseteq \mathbf{R}^d$ via their exponent vectors, let $P(\mathfrak{a}) \subseteq \mathbf{R}^d$ denote the convex hull of the points determined by monomials in $\mathfrak{a}$. Thus $P(\mathfrak{a})$ is a closed unbounded convex region lying in the first orthant. Howald's theorem states that for any $c > 0$, $\mathcal{J}(\mathfrak{a}^c)$ is the monomial ideal spanned by all monomials $t^v$ whose exponent vectors $v$ satisfy
$$v + \mathbf{1} \in \mathrm{interior}\left(c \cdot P(\mathfrak{a})\right),$$
where $\mathbf{1} = (1, 1, \ldots, 1)$.

Now $P(\mathfrak{a})$ is cut out in the first orthant $\{\zeta_1 \geq 0, \ldots, \zeta_d \geq 0\}$ by finitely many inequalities of the form
$$\ell_\alpha(\zeta_1, \ldots, \zeta_d) \geq 1,$$
where the $\ell_\alpha$ are linear forms with non-negative rational coefficients. Then the jumping numbers of $\mathfrak{a}$ on $\mathbf{C}^d$ (or equivalently at the origin) consist of all the rational numbers
$$\xi_v = \min_\alpha \left\{ \ell_\alpha(v + \mathbf{1}) \right\}$$
as $v$ ranges over $\mathbf{N}^d$. In fact, it follows from Howald's theorem that $\xi_v$ is the smallest rational number $c > 0$ such that $t^v \notin \mathcal{J}(c \cdot \mathfrak{a})$. Note however that different $v \in \mathbf{N}^d$ might give rise to the same coefficient $\xi$. □

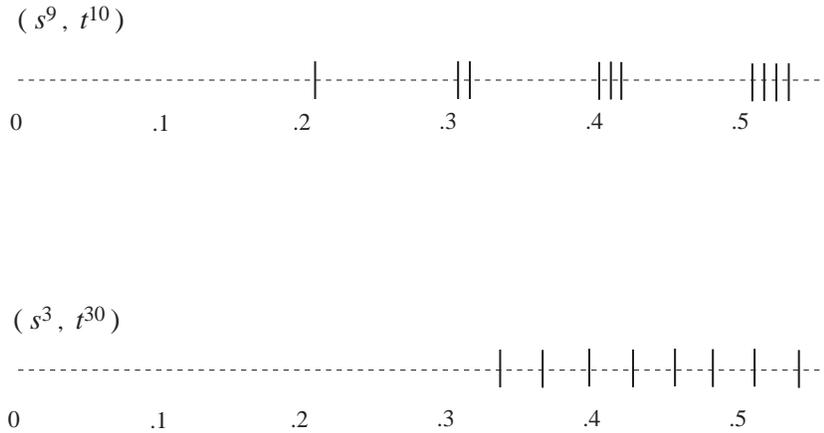

FIGURE 1. Jumping coefficients of $(s^9, t^{10})$ and $(s^3, t^{30})$.

**Example 1.9** (Diagonal ideals). As a special case of the preceding example, fix positive integers $m_1, \ldots, m_d$ and consider the ideal $\mathfrak{a} = (t_1^{m_1}, \ldots, t_d^{m_d})$. Then the non-zero jumping coefficients of $\mathfrak{a}$ on $\mathbf{C}^d$ (or at 0) consist precisely of the rational numbers
$$\frac{e_1 + 1}{m_1} + \ldots + \frac{e_d + 1}{m_d}$$
as $e_1, \ldots, e_d \in \mathbf{N}$ range over all non-negative integers. □

**Remark 1.10.** It is natural to represent jumping numbers graphically. Figure 1 shows the first few non-zero jumping numbers for the ideals $(s^9, t^{10})$ and $(s^3, t^{30})$ in $\mathbf{C}[s,t]$: the exponents are chosen so that both ideals have the same Samuel multiplicity, and so that the pictured jumping coefficients appear without repetitions in the sense of Definition 1.22. □

**Example 1.11** (Non-degenerate polynomials). Given a polynomial
$$f \in \mathbf{C}[t_1, \ldots, t_d] = \mathbf{C}[t],$$
consider the jumping numbers $\xi_i(f)$ of the principal ideal $(f)$ on $X = \mathbf{C}^d$. A second theorem of Howald [14], [13] gives an explicit condition on $f$ to guarantee that the $\xi_i(f)$ can be computed from the *term ideal* of $f$, that is, the monomial ideal $\mathfrak{a}_f \subseteq \mathbf{C}[t]$ generated by the monomials occurring in $f$. Specifically, write $P(f) = P(\mathfrak{a}_f)$ for the Newton polyhedron associated to $\mathfrak{a}_f$. Given any face $\sigma$ of $P(f)$ — including $\sigma = P(f)$ — denote by $f_\sigma$ the sum of those terms of $f$ corresponding to points lying on $\sigma$. One says that $f$ is *non-degenerate along* $\sigma$ if the 1-form $df_\sigma$ is nowhere vanishing along the torus $(\mathbf{C}^*)^d$, and one says that $f$ is *non-degenerate* if it is non-degenerate along each face of $P(f)$.[4] Howald proves that if $f$ is non-degenerate, then
$$\mathcal{J}(c \cdot f) = \mathcal{J}(c \cdot \mathfrak{a}_f) \tag{3}$$

---

[4] If $f$ has a non-zero constant term, then $P(f)$ contains the origin $\{0\}$ as a face. In this case we agree that $f$ is non-degenerate along $\{0\}$ even though $f_{\{0\}}$ is constant and $df_{\{0\}} = 0$.

for every $0 < c < 1$. Observing that $\mathcal{J}(c \cdot f) = (f) \cdot \mathcal{J}((c-1) \cdot f)$ when $c \geq 1$, it follows that if $f$ is non-degenerate then its jumping numbers on $\mathbf{C}^d$ are computed by the algorithm in 1.8. Howald also gives an analogous condition under which (3) holds in a neighborhood of the origin. $\square$

We now establish some formal properties of jumping coefficients. As above, $X$ is a smooth complex variety of dimension $d$.

**Proposition 1.12** (Periodicity)**.** *Let $\mathfrak{a} \subseteq \mathcal{O}_X$ be a non-trivial ideal, and let $x$ be any fixed point in the closed subscheme of $X$ defined by $\mathfrak{a}$. If $\xi > d - 1$ then $\xi$ is a jumping number for $\mathfrak{a}$ at $x$ if and only if $(\xi + 1)$ is a jumping number for $\mathfrak{a}$ at $x$.*

**Example 1.13.** Let $\mathfrak{m} \subseteq \mathcal{O}_X$ be the maximal ideal of $X$ at $x$. Then the jumping numbers of $\mathfrak{m}$ consist of all positive integers $\geq d$. This shows that one cannot remove the restriction $\xi > d - 1$ in the Proposition. $\square$

*Proof of Proposition 1.12.* The question is local, so we may assume that $X$ is an affine variety. Suppose that $\xi$ is a jumping number of $\mathfrak{a}$ at $x$. Then we can find a regular function $h$ of $X$ such that $h \in \mathcal{J}(\mathfrak{a}^{\xi-\epsilon})$ for every $\epsilon > 0$, but such that $h$ is not in $\mathcal{J}(\mathfrak{a}^\xi)$. Let $g \in \mathfrak{a}$ be a general element. Using the definition one checks that $g \cdot h \in \mathcal{J}(\mathfrak{a}^{\xi+1-\epsilon})$ for every $\epsilon > 0$, but that $g \cdot h$ is not in $\mathcal{J}(\mathfrak{a}^{\xi+1})$. This shows that $\xi + 1$ is also a jumping number for $\mathfrak{a}$.

For the converse, fix an integer $m \geq d - 1$. A theorem of Skoda (cf. [4, Theorem 11.17] or [20, Chapter 12]) asserts that
$$\mathcal{J}(\mathfrak{a}^{m+\xi+1}) = \mathfrak{a} \cdot \mathcal{J}(\mathfrak{a}^{m+\xi})$$
for any rational $\xi > 0$. This shows that if $m + \xi$ is not a jumping number of $\mathfrak{a}$, then $m + \xi + 1$ is also not a jumping number of $\mathfrak{a}$. $\square$

**Remark 1.14.** The argument just completed shows that if $\xi$ is a jumping coefficient for $\mathfrak{a}$ (possibly with $\xi \leq d - 1$) then $\xi + 1$ is also a jumping number. $\square$

**Remark 1.15** (Periodicity for integral divisors)**.** If $A$ is an integral divisor on $X$, then $\mathcal{J}((\xi+1) \cdot A) = \mathcal{J}(\xi \cdot A) \otimes \mathcal{O}_X(-A)$ for any $\xi > 0$. It follows as above that the collection of jumping numbers of an integral divisor is periodic with period 1. $\square$

**Remark 1.16.** One can interpolate between the statements in 1.12 and 1.15 by taking into account the number of generators of $\mathfrak{a}$ in a neighborhood of $x$. In fact, suppose that $X$ is affine and that $\mathfrak{a}$ — or more generally a reduction of $\mathfrak{a}$ — is generated by $p$ functions. Then a variant of Skoda's theorem states that $\mathcal{J}(\mathfrak{a}^{m+c+1}) = \mathfrak{a} \cdot \mathcal{J}(\mathfrak{a}^{m+c})$ for $c > 0$ and $m \geq p - 1$ (cf. [20, Chapter 12]). So in this case we get the periodicity of jumping numbers $\xi > p - 1$. $\square$

The subadditivity theorem of [5] leads to an inequality among jumping coefficients:

**Proposition 1.17.** *Let $\mathfrak{a} \subseteq \mathcal{O}_X$ be a non-trivial ideal, and $x \in X$ a point in the zeroes of $\mathfrak{a}$. Then for every $i \geq 1$ the jumping numbers of $\mathfrak{a}$ at $x$ satisfy*
$$\xi_{i+1}(\mathfrak{a}; x) \leq \xi_1(\mathfrak{a}; x) + \xi_i(\mathfrak{a}; x).$$

*Proof.* Write $\xi_j = \xi_j(\mathfrak{a}; x)$. Then $\mathcal{J}(\mathfrak{a}^{\xi_1 + \xi_i}) \subseteq \mathcal{J}(\mathfrak{a}^{\xi_1}) \cdot \mathcal{J}(\mathfrak{a}^{\xi_i})$ thanks to the subadditivity theorem of [5]. In particular,
$$\mathcal{J}(\mathfrak{a}^{\xi_1 + \xi_i})_x \subseteq \mathcal{J}(\mathfrak{a}^{\xi_1})_x \cdot \mathcal{J}(\mathfrak{a}^{\xi_i})_x \subsetneq \mathcal{J}(\mathfrak{a}^{\xi_i}),$$
and it follows that $\xi_{i+1} \le \xi_1 + \xi_i$. □

**Remark 1.18.** It would be interesting to know whether jumping numbers satisfy any other universal inequalities. □

**Remark 1.19.** One can use jumping numbers to bound some of the data associated to a resolution of singularities. Specifically, let $\mu : X' \longrightarrow X$ be a log resolution of $\mathfrak{a}$, so that $\mathfrak{a} \cdot \mathcal{O}_{X'} = \mathcal{O}_{X'}(-F)$ for $F = \sum a_i E_i$ an effective divisor with simple normal crossing support. Note that if $\xi$ is a jumping number for $\mathfrak{a}$, then $\xi \cdot a_i$ is an integer for some $i$. Let $m$ be the least common multiple of the $a_i$'s. Then clearly $\xi_{i+1} - \xi_i \ge \frac{1}{m}$. Thus if two consecutive jumping numbers are close, then the integer $m$ must be large. □

We next use a theorem of Mustaţă [28] to prove a statement of Thom-Sebastiani type for a sum of ideals in independent sets of variables.

**Proposition 1.20.** *Let $X$ and $Y$ be smooth varieties, and fix points $x \in X$, $y \in Y$. Given ideal sheaves $\mathfrak{a} \subseteq \mathcal{O}_X$, $\mathfrak{b} \subseteq \mathcal{O}_Y$ on $X$ and $Y$ respectively, denote by $(\mathfrak{a}, \mathfrak{b})$ the ideal they generate on $X \times Y$, i.e.*
$$(\mathfrak{a}, \mathfrak{b}) =_{def} p^{-1}\mathfrak{a} + q^{-1}\mathfrak{b} \subseteq \mathcal{O}_{X \times Y}$$
*($p$ and $q$ being the projections of $X \times Y$ onto its factors). Then the non-zero jumping numbers of $(\mathfrak{a}, \mathfrak{b})$ at $(x, y)$ consist precisely of the sums*
$$(4) \qquad \xi_i(\mathfrak{a}; x) + \xi_j(\mathfrak{b}; y) \qquad (i, j \ge 1)$$
*of the jumping numbers of $\mathfrak{a}$ and $\mathfrak{b}$ at $x$ and $y$ respectively.*

**Example 1.21.** The jumping coefficients of the ideal $(t^m) \subseteq \mathbf{C}[t]$ in one variable consist of all rational numbers $\{\frac{e+1}{m}\}_{e \ge 0}$. Hence one recovers from the Proposition the computation of the jumping coefficients of the "diagonal" ideal $(t_1^{m_1}, \ldots, t_d^{m_d})$ appearing in Example 1.9. (Compare [6, Example 2.8].) □

*Proof of Proposition 1.20.* Mustaţă [28, Theorem 0.3] proves that for any rational $c > 0$:
$$(5) \qquad \mathcal{J}(X \times Y, (\mathfrak{a}, \mathfrak{b})^c) = \sum_{\lambda + \mu = c} p^{-1} \mathcal{J}(X, \mathfrak{a}^\lambda) \cdot q^{-1} \mathcal{J}(Y, \mathfrak{b}^\mu).$$
It follows right away any non-zero jumping coefficient of $(\mathfrak{a}, \mathfrak{b})$ appears among the quantities occurring in (4). To go the other way, let $\xi, \xi'$ be jumping numbers of $\mathfrak{a}, \mathfrak{b}$ at $x$ and $y$ respectively, and set $\delta = \xi + \xi'$: we need to show that for all $\varepsilon > 0$
$$(*) \qquad \mathcal{J}(X \times Y, (\mathfrak{a}, \mathfrak{b})^\delta) \subsetneq \mathcal{J}(X \times Y, (\mathfrak{a}, \mathfrak{b})^{\delta - 2\varepsilon})$$
in a neighborhood of $(x, y)$. For this we use another computation of Mustaţă's [28, Lemma 1.2], namely that
$$(6) \qquad \sum_{\lambda + \mu = c} p^{-1} \mathcal{J}(X, \mathfrak{a}^\lambda) \cdot q^{-1} \mathcal{J}(Y, \mathfrak{b}^\mu) = \bigcap_{\lambda + \mu = c} \left( \mathcal{J}(X, \mathfrak{a}^\lambda), \mathcal{J}(Y, \mathfrak{b}^\mu) \right)$$

for any $c > 0$. Assuming as we may that $X$ and $Y$ are affine, take any functions
$$f \in \mathcal{J}(\mathfrak{a}^{\xi-\varepsilon}) \subseteq \mathbf{C}[X] \ , \ g \in \mathcal{J}(\mathfrak{b}^{\xi'-\varepsilon}) \subseteq \mathbf{C}[Y].$$
Then $p^*f \cdot q^*g \in \mathcal{J}((\mathfrak{a},\mathfrak{b})^{\delta-2\varepsilon})$ thanks to (5). On the other hand, using the isomorphism
$$\frac{\mathbf{C}[X \times Y]}{(\mathfrak{j},\mathfrak{j}')} = \frac{\mathbf{C}[X]}{\mathfrak{j}} \otimes_{\mathbf{C}} \frac{\mathbf{C}[Y]}{\mathfrak{j}'}$$
for any ideals $\mathfrak{j} \subseteq \mathbf{C}[X]$ and $\mathfrak{j}' \subseteq \mathbf{C}[Y]$, we see that
$$p^*f \cdot q^*g \in \left( \mathcal{J}(\mathfrak{a}^\xi), \mathcal{J}(\mathfrak{b}^{\xi'}) \right) \quad \text{if and only if} \quad f \in \mathcal{J}(\mathfrak{a}^\xi) \text{ or } g \in \mathcal{J}(\mathfrak{b}^{\xi'}).$$
So if we take $f \notin \mathcal{J}(\mathfrak{a}^\xi)$ and $g \notin \mathcal{J}(\mathfrak{b}^{\xi'})$ then it follows from (6) that $p^*f \cdot q^*g \notin \mathcal{J}((\mathfrak{a},\mathfrak{b})^\delta)$, giving (*). □

Finally, we discuss the situation for ideals of finite colength. Assume henceforth that $\mathfrak{a} \subseteq \mathcal{O}_X$ vanishes at a single point $x \in X$, and denote by $\mathfrak{m} \subseteq \mathcal{O}_X$ the maximal ideal of $x$. Then all of the multiplier ideals $\mathcal{J}(\mathfrak{a}^c)$ are $\mathfrak{m}$-primary, and hence their stalks at $x$ have finite codimension (as a $\mathbf{C}$-vector space) in the local ring $\mathcal{O}_x X = \mathcal{O}_{X,x}$. This allows one to assign in the natural way a multiplicity to each of the jumping numbers $\xi_i = \xi_i(\mathfrak{a};x)$:

**Definition 1.22** (Multiplicity of a jumping number). Still assuming that $\mathfrak{a}$ is $\mathfrak{m}$-primary, the multiplicity attached to the jumping number $\xi_i = \xi_i(\mathfrak{a};x)$ is the codimension of $\mathcal{J}(\mathfrak{a}^{\xi_i})_x$ in $\mathcal{J}(\mathfrak{a}^{\xi_{i-1}})_x$. We denote by
$$(7) \qquad 0 = \kappa_0(\mathfrak{a}) \leq \kappa_1(\mathfrak{a}) \leq \kappa_2(\mathfrak{a}) \leq \ldots$$
the jumping numbers of $\mathfrak{a}$ at $x$, each repeated according to its multiplicity. □

Note that if $\mathcal{J}(\mathfrak{a}^{\kappa_\ell}) \supsetneq \mathcal{J}(\mathfrak{a}^{\kappa_{\ell+1}})$, then $\mathcal{J}(\mathfrak{a}^{\kappa_\ell})_x$ has colength exactly $\ell$ in $\mathcal{O}_x X$.

**Example 1.23.** Consider the maximal ideal $\mathfrak{a} = (s,t) \subseteq \mathbf{C}[s,t]$ in the polynomial ring in two variables. Here one has
$$\kappa_1 = 2, \quad \kappa_2 = \kappa_3 = 3, \quad \kappa_4 = \kappa_5 = \kappa_6 = 4, \quad \text{etc.} \qquad □$$

The coefficients $\kappa_i$ satisfy a semicontinuity property analogous to that obeyed by the spectrum of a singularity [32].

**Proposition 1.24** (Semicontinuity). *Let $T$ be a smooth curve, and let $\mathfrak{a} \subseteq \mathcal{O}_{X \times T}$ be an ideal whose zeroes are supported on $\{x\} \times T$ for some $x \in X$. Given $t \in T$ write*
$$\mathfrak{a}_t \subseteq \mathcal{O}_{X \times \{t\}} = \mathcal{O}_X$$
*for the ideal obtained as the fibre of $\mathfrak{a}$ over $t$. Then there is a dense Zariski-open set $U \subseteq T$ such that all of the jumping numbers $\kappa_i(\mathfrak{a}_t)$ are independent of $t$ provided that $t \in U$. Moreover if $t \in U$ and $t^* \in T$ is an arbitrary point, then*
$$(8) \qquad \kappa_i(\mathfrak{a}_{t^*}) \leq \kappa_i(\mathfrak{a}_t) \quad \text{for all} \quad i \geq 0.$$

In other words, the jumping coefficients $\kappa_i(\mathfrak{a}_t)$ are generically constant, and can only decrease as a general point $t \in T$ specializes to $t^* \in T$.

*Proof.* Given a positive number $\lambda \in \mathbf{Q}$, put
$$K(\lambda, t) \;=\; \max\{j \mid \kappa_j(\mathfrak{a}_t) \le \lambda\}.$$

The first statement is equivalent to the constancy of $K(\lambda, t)$ for $t \in U$ and all $\lambda \ge 0$, while (8) is equivalent to the assertion that

(9) $$K(\lambda, t^*) \;\ge\; K(\lambda, t) \quad \text{for all} \;\; \lambda \ge 0.$$

Bearing in mind the comment immediately following Definition 1.22, it is enough to prove that $\operatorname{codim}_{\mathcal{O}_x X}\big(\mathcal{J}(\lambda \cdot \mathfrak{a}_t)\big)$ is generically constant, whereas

(10) $$\operatorname{codim}_{\mathcal{O}_x X}\Big(\mathcal{J}(\lambda \cdot \mathfrak{a}_{t^*})\Big) \;\ge\; \operatorname{codim}_{\mathcal{O}_x X}\Big(\mathcal{J}(\lambda \cdot \mathfrak{a}_t)\Big).$$

For this consider the multiplier ideals $\mathcal{J}(X \times T, \lambda \cdot \mathfrak{a})$ on $X \times T$. There is an open subset $U_1 \subseteq T$ such that each of these ideals is flat over $U_1$.[5] Moreover by the theorem on generic restrictions (cf. [20, Chapter 9.5]) there is an open set $U_2 \subseteq T$ such that

(11) $$\mathcal{J}(X \times T, \lambda \cdot \mathfrak{a}) \cdot \mathcal{O}_{X \times \{t\}} \;=\; \mathcal{J}(X, \lambda \cdot \mathfrak{a}_t)$$

for general $t \in U_2$ and every $\lambda > 0$. Then $\operatorname{codim}_{\mathcal{O}_x X}\big(\mathcal{J}(\lambda \cdot \mathfrak{a}_t)\big)$ is constant for $t \in U_1 \cap U_2$. Now consider an arbitrary point $t^* \in T$. One has
$$\mathcal{J}(X, \lambda \cdot \mathfrak{a}_{t^*}) \;\subseteq\; \mathcal{J}(X \times T, \lambda \cdot \mathfrak{a}) \cdot \mathcal{O}_{X \times \{t^*\}}$$
for every $\lambda > 0$ by the Restriction Theorem (cf. [5, §1] or [20, Chapter 9.5]). On the other hand, the codimension in $\mathcal{O}_X = \mathcal{O}_{X \times \{t\}}$ of the ideals $\mathcal{J}(X \times T, \lambda \cdot \mathfrak{a}) \cdot \mathcal{O}_{X \times \{t\}}$ is in any event upper semicontinuous in $t$ and (10) follows.  □

**Remark 1.25.** Proposition 1.24 could be deduced from results of Steenbrink [32] via the connection with the spectrum of a singularity (Remark 2.2). However we felt it worthwhile to give a direct elementary proof. Furthermore, in §5 the same argument will give semicontinuity for the jumping numbers attached to a graded system of finite colength ideals (Proposition 5.14).  □

**Remark 1.26** (Analytic ideals)**.** Viewing $X$ as a complex manifold, consider any coherent analytic ideal sheaf $\mathfrak{a} \subseteq \mathcal{O}_X$. The multiplier ideals $\mathcal{J}(X, \mathfrak{a}^c) \subseteq \mathcal{O}_X$ can be defined as above (or else via a direct analytic construction as in Remark 1.2), and Lemma 1.3 goes through with no change. So given any $x \in X$ the local jumping coefficients $\xi_i(\mathfrak{a}; x)$ are defined as before; in particular, a germ at $x$ of an analytic function $f$ gives rise local jumping numbers $\xi_i(f; x)$. Propositions 1.12, 1.17, 1.20 and 1.24 likewise remain valid in this analytic setting except that in the extension of 1.24, $T$ should be taken to be a small disk and $U$ should be taken to be the complement of a discrete subset of $T$.  □

---

[5]It is enough to take $U_1$ with the property that no exceptional divisor in a log resolution of $\mathfrak{a}$ on $X \times T$ maps to a point in $U_1$.

## 2. Bernstein-Sato Polynomials

In this section we establish Theorem B from the Introduction.

We start by recalling the set-up. Fix a non-zero polynomial $f \in \mathbf{C}[t_1, \ldots, t_d] = \mathbf{C}[t]$ on $\mathbf{C}^d$, and let $s$ be a new variable. Then there exists a differential operator $P$ whose coefficients are polynomials in $s$ and $t_1, \ldots, t_d$, together with a non-zero polynomial $b(s) \in \mathbf{C}[s]$, satisfying the formal identity:

$$(12) \qquad Pf^{s+1} \;=\; b(s) \cdot f^s.$$

The set of all polynomials $b(s) \in \mathbf{C}[s]$ for which such an identity holds (for some operator $P$) forms an ideal, and the unique monic generator for this ideal is called the *Bernstein-Sato polynomial* of $f$, sometimes written $b_f(s)$. We refer to [26], [17], [37], [18] for proofs, references and further information. The situation remains the same if one takes $f$ to be the germ of a holomorphic function at a point $x \in \mathbf{C}^d$, except that now the coefficients of $P$ are polynomial in $s$ and analytic in $t$. Note that one can interpret (12) as an equality of distributions provided that both sides make sense as distributions.

It was established by Yano [37], Lichten [22] and Kollár [18] that if $\operatorname{lct}(f)$ is the log-canonical threshold of $f$ on $\mathbf{C}^d$, then $-\operatorname{lct}(f)$ is the largest root $b_f(s)$. The result for which we are aiming generalizes this fact:

**Theorem 2.1.** *Let $\xi$ be a jumping coefficient of $f$ on $\mathbf{C}^d$ which lies in the interval $(0,1]$. Then $b_f(-\xi) = 0$.*

The analogous statement holds when $f$ is the germ of a holomorphic function at some $x \in \mathbf{C}^d$: here $\xi$ should be a local jumping number of $f$ at $x$. We recall that in any event $b_f(-1) = 0$.

**Remark 2.2** (Relation to the spectrum of a singularity)**.** At least when $f$ has isolated singularities, Loeser and Vaquié [25] observe that results of Varchenko [36] establish a connection between jumping numbers and the *spectrum* of $f$: this in turn implies the Theorem in the case of isolated singularities. Specifically, Steenbrink and Varchenko constructed a mixed Hodge structure on the $(d-1)$-st cohomology group of the Milnor fiber of $f$. They then defined $\operatorname{Sp}(f)$, the spectrum of $f$, in the following manner. For $p$ a fixed integer between $0$ and $d-1$, let $\alpha$ be a rational number such that $d - p - 2 < \alpha \leq d - p - 1$. Then $\alpha \in \operatorname{Sp}(f)$ if and only if $\exp(2\pi i \alpha)$ is an eigenvalue of the monodromy operator on $\frac{F^p H}{F^{p+1} H}$, where $F^p H$ is the $p$-th piece of the Hodge filtration on the $(d-1)$-th cohomology group of the Milnor fiber of $f$. Vaquié [33, p. 1191], [34, p. 390] explicitly observed that the result of Varchenko implies that for $\alpha$ between $-1$ and $0$, $\alpha \in \operatorname{Sp}(f)$ if and only if $\alpha + 1$ is a jumping number of $f$. The case of non-isolated singularities is treated by Budur [2]. □

**Remark 2.3.** It is likely that jumping coefficients are related to other invariants of singularities that appear in the literature, for example those studied by Saito in [30]. It would be interesting to know the precise connections. □

In general, not all roots of $b_f(s)$ lying in $[-1, 0)$ come from jumping numbers (Example 2.5). However combined with Proposition 1.17, Theorem 2.1 gives:

**Corollary 2.4.** *With $f$ as above, let*
$$0 > r_1 > r_2 > \ldots > r_t = -1$$
*be the distinct roots of $b_f(s)$ lying between $0$ and $-1$. Then*

(13) $$r_{i+1} \geq r_i + r_1$$

*for every $1 \leq i < t$.*

Again this holds also when $f$ is a holomorphic germ. As we noted in the Introduction, van Doorn and Steenbrink [35] have established some related inequalities in the case of isolated singularities.

*Proof of Corollary 2.4.* Put $\rho_i = -r_i$, so that the question is to establish the inequality $\rho_{i+1} \leq \rho_i + \rho_1$ provided that $i < t$. Let $\xi'$ be the greatest jumping number of $f$ which is $\leq \rho_i$, and let $\xi$ be the next jumping coefficient after $\xi'$, so that $\xi' \leq \rho_i < \xi$. One has $\rho_1 = \xi_1(f) = \text{lct}(f)$ thanks the theorem of Yano-Lichten-Kollár, and so Proposition 1.17 gives
$$\xi \leq \xi' + \rho_1 \leq \rho_i + \rho_1.$$
Now $\xi' \leq \rho_i < 1$ and hence $\xi \leq 1$ (since in any event 1 is a jumping coefficient of $f$). Therefore Theorem 2.1 applies to show that $\xi = \rho_j$ for some $j > i$. Then $\rho_{i+1} \leq \rho_j$, and (13) follows. □

**Example 2.5** (Non-jumping roots). Saito [30, Example 3.5, p. 69] shows that if $f$ is the polynomial $f = x^5 + y^4 + x^3 y^2$, then $b_f(s)$ has roots in $[-1, 0)$ which do not come from the spectrum of $f$. In view of Remark 2.2, this gives an example where not every root of the Bernstein polynomial arises from a jumping number. □

We now turn to the demonstration of Theorem 2.1, closely following Kollár's proof in [18, §10] of the theorem of Yano-Lichten-Kollár. For the benefit of the reader not versed in such matters, we go through the argument in some detail.

*Proof of Theorem 2.1.* Fix any jumping number $\xi \in (0, 1]$ of $f$ on $\mathbf{C}^d$ and let $\xi' < \xi$ be the previous jumping number (so in particular $\xi' = 0$ in case $\xi = \xi_1(f)$ is the log canonical threshold). Thus

(14) $$\mathcal{J}(c \cdot f) = \mathcal{J}(\xi' \cdot f)$$

for all $c \in [\xi', \xi)$, but

(15) $$\mathcal{J}(\xi \cdot f)_x \subsetneq \mathcal{J}(\xi' \cdot f)_x$$

for some point $x \in \mathbf{C}^d$ at which $f$ vanishes. In particular, given $c \in [\xi', \xi)$ and a germ $g \in \mathcal{J}(c \cdot f)_x$, it follows from (14) and Remark 1.2 that there exists a small ball $B$ around $x$ — which by the coherence of multiplier ideals one can take independent of $c$ – such that

(16) $$\int_B \frac{|g|^2}{|f^2|^c} < \infty.$$

On the other hand, if $g \notin \mathcal{J}(\xi \cdot f)_x$ then $\frac{|g|^2}{|f^2|^\xi}$ is not integrable at $x$. The plan roughly speaking is to exploit the fact that the integrals in (16) must consequently become unbounded as $c$ approaches $\xi$.

Turning to details, write $b(s) = b_f(s)$ and fix $P$ as in (12). As explained in [, p. ], one has
$$(Pf^{s+1}) \cdot (\bar{P}\bar{f}^{s+1}) = P\bar{P}(|f^2|^{s+1})$$
thanks to the fact the holomorphic and anti-holomorphic operators commute. This gives the formal identity
$$P\bar{P}(|f^2|^{s+1}) = |b(s)|^2 \cdot |f^2|^s,$$
leading to the relation
$$(17) \quad |g|^2 \cdot P\bar{P}(|f^2|^{s+1}) = |g|^2 \cdot |b(s)|^2 \cdot |f^2|^s$$
for any $g$. Now take
$$s = -c \quad \text{for} \quad c \in [\xi', \xi) \qquad \text{and} \qquad g \in \mathcal{J}(c \cdot f)_x = \mathcal{J}(\xi' \cdot f)_x.$$
We claim that then both sides of (17) determine well-defined distributions on the ball $B$ appearing in (16), so that (17) holds as an equality of distributions. In fact $|g|^2 |b(s)|^2 |f^2|^s$ is integrable on $B$, and so the right side of (17) is well-defined as a distribution; and since $s + 1 \geq 0$ the same is true of the left-hand side. In particular, for any positive smooth compactly supported test function $\phi$ on $B$ we have
$$(18) \quad \int_B \phi \cdot |g^2| \cdot |b(-c)|^2 \cdot |f^2|^{-c} = \int_B |f^2|^{1-c} \cdot P\bar{P}(\phi|g^2|).$$
One sees that the right hand side of (18) is uniformly bounded above for all $c \in [\xi', \xi)$ by a finite positive number $M$ depending on $\phi$ and $g$. On the other hand, if $\phi$ takes the constant value 1 on some ball $B' \subset B$ about $x$ then the left-hand side of (18) is at least
$$\int_{B'} \frac{|g^2||b(-c)|^2 \phi}{|f^2|^c} = |b(-c)|^2 \int_{B'} \frac{|g^2|}{|f^2|^c}.$$
All told, we have
$$(19) \quad |b(-c)|^2 \cdot \int_{B'} \frac{|g^2|}{|f^2|^c} \leq M < \infty$$
for every $c \in [\xi', \xi)$.

Now since $\xi > \xi'$ is a jumping number for $f$, we can take
$$g \in \mathcal{J}(\xi' \cdot f)_x \qquad \text{with} \qquad g \notin \mathcal{J}(\xi \cdot f)_x.$$
Thus $\frac{|g|^2}{|f^2|^\xi}$ fails to be integrable at $x$, and consequently the integrals $\int_{B'} \frac{|g^2|}{|f^2|^c}$ must go to infinity as $c$ approaches $\xi$ from below thanks to the monotone convergence theorem. But then it follows from (19) that $b(-\xi) = 0$, as required. $\square$

3. Uniform Artin-Rees Numbers

In this section we establish a connection between jumping coefficients and Huneke's uniform Artin-Rees theorem [15].

Let $X$ be a smooth variety of dimension $d$, and let $\mathfrak{a} \subseteq \mathcal{O}_X$ be an ideal sheaf on $X$. Our main technical result is the following:

**Theorem 3.1.** *Let $\xi' < \xi$ be consecutive jumping coefficients of $\mathfrak{a}$ on $X$, i.e. assume that $\xi' = \xi_i(\mathfrak{a}; X)$ and $\xi = \xi_{i+1}(\mathfrak{a}; X)$ for some index $i$. Then*
$$\text{(20)} \qquad \mathfrak{b}^m \cdot \mathcal{J}(\mathfrak{a}^{\xi'}) \cap \mathcal{J}(\mathfrak{a}^\xi) \ \subseteq \ \mathfrak{b}^{m-d} \cdot \mathcal{J}(\mathfrak{a}^\xi)$$
*for any ideal sheaf $\mathfrak{b} \subseteq \mathcal{O}_X$ and all $m \geq d$.*

In the terminology of [15], the Theorem asserts that the dimension $d$ of $X$ is a uniform Artin-Rees number for a pair $\mathcal{J}(\mathfrak{a}^\xi) \subseteq \mathcal{J}(\mathfrak{a}^{\xi'})$ of consecutive multiplier ideals. Given a fixed point $x \in X$, the analogue of (20) holds for stalks at $x$ when $\xi'$ and $\xi$ are consecutive jumping numbers at $x$. The proof appears at the end of the section.

By considering the neighboring pairs of ideals in the chain
$$\mathcal{O}_X \ \supseteq \ \mathcal{J}(\mathfrak{a}^{\xi_1}) \ \supseteq \ \mathcal{J}(\mathfrak{a}^{\xi_2}) \ \supseteq \ \ldots \ \supseteq \mathcal{J}(\mathfrak{a}^{\xi_\ell}),$$
one then arrives at an effective uniform Artin-Rees number for any of the multiplier ideals $\mathcal{J}(\mathfrak{a}^c)$:

**Corollary 3.2.** *Given any rational number $c > 0$, let $\ell = \ell(\mathfrak{a}, c; X)$ be the number of non-zero jumping coefficients of $\mathfrak{a}$ on $X$ that are $\leq c$. Then*
$$\text{(21)} \qquad \mathfrak{b}^m \cap \mathcal{J}(\mathfrak{a}^c) \ \subseteq \ \mathfrak{b}^{m-\ell d} \cdot \mathcal{J}(\mathfrak{a}^c)$$
*for every ideal $\mathfrak{b} \subseteq \mathcal{O}_X$ and every $m \geq \ell d$.* □

The integer $\ell = \ell(\mathfrak{a}, c; X)$ is thus the length of the maximal length chain of the form
$$\mathcal{O}_X \ \supsetneq \ \mathcal{J}(\mathfrak{a}^{\xi_1}) \ \supsetneq \ \mathcal{J}(\mathfrak{a}^{\xi_2}) \ \supsetneq \ \ldots \ \supsetneq \ \mathcal{J}(\mathfrak{a}^c).$$
We will call this the *jumping length* of the multiplier ideal $\mathcal{J}(\mathfrak{a}^c)$. Note that a given ideal $\mathfrak{j}$ may be a multiplier ideal in more than one way: it is possible that $\mathfrak{j} = \mathcal{J}(\mathfrak{a}^c) = \mathcal{J}(\mathfrak{b}^d)$ for different ideals $\mathfrak{a}, \mathfrak{b}$ and positive rational numbers $c, d$. In this case, if we do not refer explicitly to $\mathfrak{j}$ as $\mathcal{J}(\mathfrak{a}^c)$, then we interpret the jumping length of the ideal $\mathfrak{j}$ to be the minimal possible jumping length for any way of representing the ideal $\mathfrak{j}$ as a multiplier ideal.

Again, the analogous statement holds for stalks and jumping coefficients at a fixed point $x \in X$. Stated in local commutative algebra language, therefore, we have the following corollary.

**Corollary 3.3.** *Let $\mathfrak{j}$ be any multiplier ideal in a d-dimensional regular local ring $R$ essentially of finite type over a field of characteristic zero, and let $\ell$ be its jumping length. Then for any ideal $\mathfrak{b}$ of $R$, we have*
$$\mathfrak{b}^m \cap \mathfrak{j} \ \subseteq \ \mathfrak{b}^{m-\ell d} \cdot \mathfrak{j}$$

*for every $m \geq d\ell$.*

**Remark 3.4** (Realizing multiplier ideals)**.** The previous corollaries raise the very interesting question: which integrally closed ideals arise as multiplier ideals? In other words, if $\mathfrak{j} \subseteq \mathcal{O}_X$ is integrally closed, when does there exist an ideal $\mathfrak{a}$ and a rational number $c > 0$ such that $\mathfrak{j} = \mathcal{J}(\mathfrak{a}^c)$? Favre-Jonsson [10] and Lipman-Watanabe [23] have recently established that in dimension 2, any such $\mathfrak{j}$ can be realized as a multiplier ideal provided that it vanishes only on a finite set. In all dimensions, every principal ideal is a multiplier ideal (Example 1.7) and as Mustaţă remarks it follows from Howald's theorem (Example 1.8) that all integrally closed monomial ideal are multiplier ideals. It seems unlikely in general that essentially every integrally closed $\mathfrak{j}$ should be a multiplier ideal, but it remains an interesting open problem to produce an example of one which isn't. $\square$

**Example 3.5** (Smooth subvarieties)**.** Let $\mathfrak{p} \subseteq \mathcal{O}_X$ be the ideal of a non-singular subvariety of $X$ having codimension $e$. Then $\mathcal{J}(\mathfrak{p}^c) = \mathfrak{p}^{[c+1-e]}$ for all $c > 0$. In particular, for any integer $r$, the ideal $\mathfrak{p}^r (= \mathcal{J}(\mathfrak{p}^{r-1+e}))$ is a multiplier ideal with jumping length $r$. The Corollary thus guarantees that

$$\mathfrak{b}^m \cap \mathfrak{p}^r \subseteq \mathfrak{b}^{m-rd} \cdot \mathfrak{p}^r$$

for every ideal $\mathfrak{b}$. When $r = 1$, Huneke [16, p. 77] established the stronger statement that in fact $d - e$ is a uniform Artin-Rees number for $\mathfrak{p}$. $\square$

The above corollaries are particularly interesting in the case of a principal ideal. Assuming for concreteness that $X$ is affine, fix a non-zero function $f \in \mathbf{C}[X]$. In this case $\mathcal{J}(f) = (f)$ and $\xi = 1$ is a jumping number for $f$ (Example 1.7). As mentioned in the Introduction, the *jumping length* $\ell(f) = \ell(f; X)$ of $f$ on $X$ in this case is thus the number of non-zero jumping coefficients of $(f)$ that are $\leq 1$.

**Example 3.6.** Consider the polynomial $f = s^3 + t^4 \in \mathbf{C}[s,t]$. Then $f$ is non-degenerate with respect to its Newton polytope, so by Howald's theorem (Example 1.11) the jumping numbers of $f$ which are less than 1 coincide with those of its term ideal $(s^3, t^4)$. Using Example 1.9 we find the first few jumping coefficients of $f$ to be:

$$0 \;<\; \tfrac{1}{3}+\tfrac{1}{4} \;<\; \tfrac{1}{3}+\tfrac{2}{4} \;<\; \tfrac{2}{3}+\tfrac{1}{4} \;<\; 1.$$

So here $\ell(f, \mathbf{C}^2) = 4$. $\square$

In this setting 3.2 yields:

**Corollary 3.7.** *Given $f \in \mathbf{C}[X]$, let $\ell = \ell(f; X)$ be the jumping length of $f$ on $X$. Then*

$$\mathfrak{b}^m \cap (f) \;\subseteq\; (f) \cdot \mathfrak{b}^{m-d\ell}$$

*for every ideal $\mathfrak{b} \subseteq \mathbf{C}[X]$ and every $m \geq d \cdot \ell$.* $\square$

Once again, the analogous local statement holds at a fixed point $x \in X$.

When $f$ has an isolated singularity, one can relate the jumping length to other invariants. Denote by $\mathrm{Jac}(f) \subseteq \mathcal{O}_X$ the Jacobian ideal of $f$, i.e. the ideal locally generated by $f$ and its partials[6]:
$$\operatorname{Jac}(f) =_{\text{locally}} \left( f, \tfrac{\partial f}{\partial x_1}, \ldots, \tfrac{\partial f}{\partial x_d} \right),$$

$x_1, \ldots, x_d$ being local coordinates of $X$. If $x \in X$ is an isolated singularity of $(f = 0)$ then $\operatorname{Jac}(f)_x$ has finite codimension in $\mathcal{O}_x X$, and the *Tyurina number* of $\tau = \tau(f; x)$ of $f$ at $x$ is defined to be
$$\tau = \dim_{\mathbf{C}} \frac{\mathcal{O}_x X}{\operatorname{Jac}(f)_x}.$$

In §4 we will prove:

**Proposition 3.8.** *Assume that the hypersurface defined by the vanishing of $f$ has an isolated singularity at $x \in X$ and otherwise is smooth. Then*
$$\operatorname{Jac}(f) \subseteq \mathcal{J}\big((f)^{1-\varepsilon}\big) \quad \text{for all} \quad \varepsilon > 0. \tag{22}$$
*In particular, the jumping length of $f$ at $x$ satisfies the bound*
$$\ell(f; x) \leq \tau(f; x) + 1.$$

The second assertion follows from (22) upon observing that the number of jumping coefficients of $(f)$ which are $< 1$ is bounded by $\operatorname{colength}_{\mathcal{O}_x X} \mathcal{J}\big((f)^{1-\varepsilon}\big)$ for small $\varepsilon > 0$.

Combining the Proposition with Corollary 3.7, we find that $k = d \cdot \big(\tau(f; x) + 1\big)$ is a uniform Artin-Rees number for $(f)$ in the case of isolated singularities. However one can do better:

**Corollary 3.9.** *In the situation of Proposition 3.8, $k = \tau(f; x) + d$ is a uniform Artin-Rees number for $(f)$.*

*Proof of Corollary 3.9.* Let $\xi' < 1$ be the largest jumping number of $f$ at $x$ which is $< 1$, and write $\mathfrak{m}$ for the maximal ideal of $x$. Then Proposition 3.8 gives the inclusion $\operatorname{Jac}(f) \subseteq \mathcal{J}\big(f^{\xi'}\big)$. Thus for any ideal $\mathfrak{b}$ and all $m \geq d$:
$$\mathfrak{b}^m \cdot \operatorname{Jac}(f) \cap (f) \subseteq \mathfrak{b}^{m-d} \cdot (f)$$
thanks to Theorem 3.1. On the other hand, since $(f = 0)$ has an isolated singularity at $x$, $\operatorname{Jac}(f)$ is $\mathfrak{m}$-primary, of colength $\tau = \tau(f; x)$ in $\mathcal{O}_x X$. As Huneke [15, §5, p. 218] observes, this implies that $\tau$ is a uniform Artin-Rees number for $\operatorname{Jac}(f) \subseteq \mathcal{O}_X$, i.e.
$$\mathfrak{b}^m \cap \operatorname{Jac}(f) \subseteq \mathfrak{b}^{m-\tau} \cdot \operatorname{Jac}(f)$$
for all ideals $\mathfrak{b}$ and $m \geq \tau$. Putting these together we find that
$$\mathfrak{b}^m \cap (f) \subseteq \mathfrak{b}^{m-\tau} \cdot \operatorname{Jac}(f) \cap (f) \subseteq \mathfrak{b}^{m-\tau-d} \cdot (f),$$
as required. □

---

[6]This differs slightly from the possibly more usual definition of $\operatorname{Jac}(f)$ as the ideal generated only by the partials of $f$. The ideal we consider here has the advantage of being intrinsically defined, and we trust that no confusion will result. In any event the inclusions we establish for $\operatorname{Jac}(f)$ imply *a fortiori* the corresponding statements for the ideal $\big(\tfrac{\partial f}{\partial x_1}, \ldots, \tfrac{\partial f}{\partial x_d}\big)$.

**Remark 3.10** (Bounds involving the Milnor number)**.** The following remark is due to Nero Budur. We saw in Proposition 3.8 that when $(f = 0)$ has an isolated singularity at $x$ then its jumping length is bounded by its Tyurina number. Under the same hypothesis, the connection between jumping numbers and the spectrum of a singularity (Remark 2.2) leads to a usually stronger bound relating the jumping length of $f$ to its Milnor number provided that $\dim X \geq 2$. Recall that the Milnor number of the hypersurface defined by $f$ at a point $x \in X$ is defined as

$$\mu \;=\; \mu(f;x) \;=\; \dim_{\mathbf{C}} \frac{\mathcal{O}_x X}{\left(\frac{\partial f}{\partial x_1}, \ldots, \frac{\partial f}{\partial x_d}\right)},$$

where $x_1, \ldots, x_d$ are local coordinates of $X$ at $x$. Now, as mentioned in Remark 2.2, a rational number $\alpha$ between $-1$ and $0$ lies in the spectrum of $f$ if and only if $\alpha + 1$ is a jumping number of $f$. On the other hand, because of a symmetry in the spectrum [19, 8.33], at most half of the distinct elements of the spectrum of $f$ are in the range $(-1, 0)$. Since the dimension of the $(d-1)$-st cohomology group of the Milnor fiber is equal to the Milnor number $\mu$, we see that there are at most $\frac{\mu}{2}$ spectral numbers in this range. Hence $\ell(f;x) \leq 1 + \frac{\mu}{2}$, and arguing as in the proof of 3.9 we find that $k = \frac{\mu}{2} + d$ is a uniform Artin-Rees number for $f$. $\square$

Finally, we turn to the proof of the Theorem.

*Proof of Theorem 3.1.* Fix an ideal $\mathfrak{b}$, and choose a rational number $\xi' \leq c < \xi$ very close to $\xi$. Then $\mathcal{J}(\mathfrak{a}^c) = \mathcal{J}(\mathfrak{a}^{\xi'})$ and $\mathfrak{b}^m \cdot \mathcal{J}(\mathfrak{a}^c) \subseteq \mathcal{J}(\mathfrak{b}^m \cdot \mathfrak{a}^c)$. We will prove that

(23) $$\mathcal{J}(\mathfrak{b}^m \cdot \mathfrak{a}^c) \cap \mathcal{J}(\mathfrak{a}^\xi) \;\subseteq\; \mathcal{J}(\mathfrak{b}^{m-1} \cdot \mathfrak{a}^\xi).$$

Then (20) follows from Skoda's theorem that

$$\mathcal{J}(\mathfrak{b}^{m-1} \cdot \mathfrak{a}^\xi) \;\subseteq\; \mathfrak{b}^{m-d} \cdot \mathcal{J}(\mathfrak{a}^\xi)$$

provided that $m \geq d$ (cf. [20, Chapter 12.1]).[7] To prove (23), let $\mu : X' \longrightarrow X$ be a common log resolution of $\mathfrak{b}$ and $\mathfrak{a}$. Write

$$\mathfrak{b} \cdot \mathcal{O}_{X'} \;=\; \mathcal{O}_{X'}(-B) \quad \text{and} \quad \mathfrak{a} \cdot \mathcal{O}_{X'} \;=\; \mathcal{O}_{X'}(-F).$$

Given an affine open set $U \subseteq X$ it follows from the definition of multiplier ideals that a function $f \in \mathcal{O}_X(U)$ lies in $\mathcal{J}(\mathfrak{b}^m \cdot \mathfrak{a}^c)(U)$ if and only if

(24) $$\mathrm{div}(\mu^* f) + K_{X'/X} - \lceil mB + cF \rceil \;\succcurlyeq\; 0$$

over $U$, and similarly $f \in \mathcal{J}(\mathfrak{a}^\xi)(U)$ if and only if

(25) $$\mathrm{div}(\mu^* f) + K_{X'/X} - \lceil \xi F \rceil \;\succcurlyeq\; 0$$

---

[7]Recall that for ideals $\mathfrak{a}_1, \mathfrak{a}_2 \subseteq \mathcal{O}_X$ and coefficients $c_1, c_2 > 0$, one defines the "mixed" multiplier ideal $\mathcal{J}(\mathfrak{a}_1^{c_1} \mathfrak{a}_2^{c_2})$ by taking a common log resolution $\mu : X' \longrightarrow X$ of $\mathfrak{a}_1, \mathfrak{a}_2$, with $\mathfrak{a}_i \cdot \mathcal{O}_{X'} = \mathcal{O}_{X'}(-F_i)$, and setting

$$\mathcal{J}(\mathfrak{a}_1^{c_1} \mathfrak{a}_2^{c_2}) \;=\; \mu_* \mathcal{O}_{X'}(K_{X'/X} - \lceil c_1 F_1 + c_2 F_2 \rceil).$$

over $U$. This being so, it is enough to show that if $c$ is sufficiently close to $\xi$ then

(26) $$\operatorname{lcm}\bigl([\xi F], [mB + cF]\bigr) \succcurlyeq [(m-1)B + \xi F].$$

For once one knows (26) it follows from (24) and (25) that

$$\operatorname{div}(\mu^* f) + K_{X'/X} - [(m-1)B + \xi F] \succcurlyeq 0$$

whenever $f \in \mathcal{J}(\mathfrak{b}^m \cdot \mathfrak{a}^c) \cap \mathcal{J}(\mathfrak{a}^\xi)$, as required.

As for (26), let $E$ be any prime divisor appearing in the support of $B + F$. By taking $c$ sufficiently close to $\xi$, we may suppose that $\operatorname{ord}_E(cF) > \operatorname{ord}_E(\xi F) - 1$. Then

$$\operatorname{ord}_E(mB + cF) \geq \operatorname{ord}_E((m-1)B + \xi F)$$

whenever $\operatorname{ord}_E(B) > 0$, whereas of course

$$\operatorname{ord}_E(\xi F) \geq \operatorname{ord}_E((m-1)B + \xi F)$$

whenever $\operatorname{ord}_E(B) = 0$. This yields (26), and completes the proof of the Theorem. □

**Remark 3.11.** It is natural to ask whether our methods give any insight into the uniform Artin-Rees property for an arbitrary inclusion of modules $N \subseteq M$ over a regular local ring $R$ in the sense of Huneke [15]. Consider a prime filtration of $M/N$. This amounts to a sequence of modules $N = N_0 \subset N_1 \subset \cdots \subset N_t = M$ such that each quotient $N_i/N_{i-1}$ is isomorphic to a cyclic module $R/P_i$, where $P_i$ is some prime ideal ideal of $R$. It is elementary to check that if $k_i$ is a uniform Artin-Rees number for the pair $P_i \subset R$, then $k = \sum k_i$ is a uniform Artin Rees number for the pair $N \subset M$. Thus in our setting if the prime ideals $P_i$ happen to be multiplier ideals, then the sum of their jumping lengths determines a the uniform Artin-Rees number for $N \subseteq M$. This re-emphasizes the interest of a special case of a question we asked earlier: is every prime ideal in a regular local ring (essentially of finite type over a field of characteristic zero) a multiplier ideal? □

## 4. Jacobian and Multiplier Ideals

In this section we establish a relation between multiplier and Jacobian ideals, which includes Proposition 3.8 as a special case. Although the general statement is not used here, we believe that it is of independent interest.

We start by fixing notation. As before, $X$ is a smooth complex variety of dimension $d$, and $\mathfrak{a} \subseteq \mathcal{O}_X$ is an ideal. Write $Z = \operatorname{Zeroes}(\mathfrak{a})$ for the subscheme defined by $\mathfrak{a}$.

**Definition 4.1.** Given $m \geq 1$, the $m^{\text{th}}$ Jacobian ideal

$$\operatorname{Jac}_m(\mathfrak{a}) \subseteq \mathcal{O}_X$$

of $\mathfrak{a}$ is the $m^{\text{th}}$ Fitting ideal of the module of differentials $\Omega^1_Z$ of $Z$ ($\Omega^1_Z$ being considered in the natural way as an $\mathcal{O}_X$-module). □

Very concretely, suppose that $\mathfrak{a}$ is locally generated by $t$ regular functions $f_1, \ldots, f_t \in \mathcal{O}_X$ in a neighborhood of $x$, and let $x_1, \ldots, x_d$ be local coordinates at $x \in X$ (so that

$dx_1, \ldots, dx_d$ are a local basis for $\Omega^1_X$ near $x$). Consider the $d \times 2t$ matrix

$$A = \begin{pmatrix} f_1 & \cdots & f_t & \frac{\partial f_1}{\partial x_1} & \cdots & \frac{\partial f_t}{\partial x_1} \\ \vdots & \vdots & \vdots & \vdots & \vdots & \vdots \\ f_1 & \cdots & f_t & \frac{\partial f_1}{\partial x_d} & \cdots & \frac{\partial f_t}{\partial x_d} \end{pmatrix}.$$

Then $A$ is a presentation matrix for $\Omega^1_Z$, and $\mathrm{Jac}_m(\mathfrak{a})$ is the ideal locally generated by by the $m \times m$ minors of $A$. In particular, when $\mathfrak{a} = (f)$ is the principal ideal defined by a single function $f$, $\mathrm{Jac}_1(\mathfrak{a}) = \mathrm{Jac}(f)$ is the Jacobian ideal considered in the previous section.

The main result of this section relates these Jacobian ideals to the multiplier ideals associated to $\mathfrak{a}$.

**Theorem 4.2.** *Let $\mathfrak{a} \subseteq \mathcal{O}_X$ be any ideal sheaf on $X$, and fix a natural number $m$.*

(i). *Assume that the multiplier ideal $\mathcal{J}(\mathfrak{a}^m)$ cuts out a scheme of codimension $\geq m+1$ in $X$. Then*
$$\mathrm{Jac}_m(\mathfrak{a}) \subseteq \mathcal{J}(\mathfrak{a}^m).$$

(ii). *Assume that the multiplier ideal $\mathcal{J}(\mathfrak{a}^m)$ cuts out a scheme of codimension exactly $m$ in $X$. Then*
$$\mathrm{Jac}_m(\mathfrak{a}) \subseteq \mathcal{J}(\mathfrak{a}^{(1-\epsilon)m}) \quad \text{for all } 0 < \epsilon \leq 1.$$

*Proof.* We start by fixing notation. Let $\mu : X' \longrightarrow X$ be a log resolution of the ideal $\mathfrak{a}$, with $\mathfrak{a} \cdot \mathcal{O}_{X'} = \mathcal{O}_{X'}(-F)$. Write

$$F = \sum_{i=1}^r a_i E_i \quad , \quad K_{X'/X} = \sum_{i=1}^r b_i E_i$$

where the $E_i$'s are distinct irreducible smooth divisors in $X'$ and $\sum E_i$ has strict normal crossings. Recall that the multiplier ideal of $\mathfrak{a}^m$ is defined as

$$\mathcal{J}(\mathfrak{a}^m) = \mu_*\left(\mathcal{O}_{X'}\left(\sum (b_i - ma_i)E_i\right)\right).$$

We need to show that if $\delta$ is one of the minors locally generating $\mathrm{Jac}_m(\mathfrak{a}^m)$, and if the hypothesis of (i) holds, then

(27) $$\mathrm{ord}_{E_i}(\mu^*\delta) \geq ma_i - b_i \quad \text{for all } i,$$

with an analogous assertion in the setting of (ii). In what follows, we focus on one of the $E_i$'s – which we call $E$ – and we put

$$a = a_i = \mathrm{ord}_E(F) \quad , \quad b = b_i = \mathrm{ord}_E(K_{X'/X}).$$

We assume also that $E$ contributes to the computation of $\mathcal{J}(\mathfrak{a}^m)$ in the sense that $ma - b > 0$, (27) being automatic in the contrary case. We set $W = \mu(E)$, and denote by $e$ the codimension of $W$ in $X$. Fix also a general point $x \in W$, and choose coordinates $\{x_1, x_2, ..., x_d\}$ for $X$ near $x$ with the property that $x_1, ..., x_e$ generate the ideal of $W$ near $x$. Finally, we suppose as above that $\mathfrak{a}$ is generated by $t$ regular functions $f_1, \ldots, f_t$ in a neighborhood of $x$.

Consider now a natural number $s \leq m$ and subsets $L, I, K \subseteq [1, d] = \{1, 2, \ldots, d\}$ with $\#L = \#I = s$ and $\#K = (m - s)$. Then the $m \times m$ minors generating $\mathrm{Jac}_m(\mathfrak{a})$ are sums and differences of terms of the form

$$\varphi_{L,I,K} = \det\left(\left(\frac{\partial f_\ell}{\partial x_i}\right)_{\substack{\ell \in L \\ i \in I}}\right) \cdot \prod_{k \in K} f_k.$$

So it suffices to show that if the hypothesis of (i) holds, then

$$\mathrm{ord}_E\left(\mu^* \varphi_{L,I,K}\right) \geq ma - b$$

with the analogous inequality

$$\mathrm{ord}_E\left(\mu^* \varphi_{L,I,K}\right) \geq ma - b - 1 \quad \left(\geq [(1-\varepsilon)ma - b]\right)$$

in the setting of (ii). Since $\mathrm{ord}_E(\mu^* f_k) \geq a$ for every $k \in K$, it is in turn enough for (i) to prove

(28) $$\mathrm{ord}_E\left(\mu^*\left(\det\left(\frac{\partial f_\ell}{\partial x_i}\right)\right)\right) \geq sa - b,$$

with the analogous statement for (ii). We will verify this by a local calculation. In order to lighten the notation, we will assume in what follows that $L = \{1, \ldots, s\}$.

The next step is to choose convenient coordinates on $X'$. To this end, let $y$ be a general point of $E \cap \mu^{-1}(x)$. Using the theorem on generic smoothness, we can find a coordinate system $\{y_1, \ldots, y_d\}$ of $X'$ near $y$ such that $y_1$ is local equation of $E$ near $y$ and such that $y_j = \mu^*(x_j)$ for $j \geq e + 1$. With these choices, we can write

(29) $$\mu^*(f_k) = y_1^a g_k \quad \text{for} \quad k \in [1, t] \quad \text{and} \quad \mu^*(x_i) = y_1 h_i \quad \text{for} \quad i \in [1, e],$$

where the $g_k$ and $h_i$ are regular functions at $y$.

With $s \leq m$ and $I \subseteq [1, d]$ as above, let $J = \{j_1, \ldots, j_{d-s}\}$ be the complement of $I$ in $[1, d]$, and write $dx_J$ for the $(d-s)$-form $dx_{j_1} \wedge \cdots \wedge dx_{j_{d-s}}$ on $X$. Consider on $X$ the $d$-form

$$\omega = df_1 \wedge \ldots \wedge df_s \wedge dx_J.$$

Then

$$\omega = \pm \det\left(\left(\frac{\partial f_\ell}{\partial x_i}\right)_{\substack{\ell \in [1,s] \\ i \in I}}\right) \cdot dx_1 \wedge \ldots \wedge dx_d.$$

Now write

$$\mu^*(\omega) = G \cdot dy_1 \wedge \cdots \wedge dy_d$$

where $G$ is a regular function at $y$. Since $K_{X'/X}$ has order $b$ along $E$ it follows that

(30) $$\mathrm{ord}_E G = \mathrm{ord}_E\left(\mu^*\left(\det\left(\frac{\partial f_\ell}{\partial x_i}\right)\right)\right) + b.$$

The crucial point now is to bound $\mathrm{ord}_E(G)$.

**Lemma 4.3.** *Keeping notation as above, one has*

$$\mathrm{ord}_E(G) \geq (sa - 1) + (e - s).$$

Granting the Lemma for the moment, we complete the proof of the Theorem. First note that if $E$ is an exceptional divisor contributing to the computation of $\mathcal{J}(\mathfrak{a}^m)$, and if $W = \mu(E)$ has codimension $e$, then certainly the ideal $\mathcal{J}(\mathfrak{a}^m)$ has height at most $e$. Thus our assumption in (i) implies that $e \geq m+1$ and hence also that $e - s \geq 1$ for any $s \leq m$. It follows from the Lemma that $\mathrm{ord}_E(G) \geq sa$, and then (30) gives the required bound (28). Since this holds for every $E$ contributing to $\mathcal{J}(\mathfrak{a}^m)$ we deduce the desired inclusion $\mathrm{Jac}_m(\mathfrak{a}) \subseteq \mathcal{J}(\mathfrak{a}^m)$. Statement (ii) is similar. $\square$

It remains only to treat the Lemma:

*Proof of Lemma 4.3.* Referring to equation (29) we find that
$$\mu^*(df_\ell) \;=\; d(y_1^a g_\ell) \;=\; a y_1^{a-1} g_\ell dy_1 + y_1^a dg_\ell$$
for $1 \leq \ell \leq s$. Similarly for $1 \leq i \leq e$, we can write
$$\mu^*(dx_i) \;=\; d(y_1 h_i) \;=\; h_i dy_1 + y_1 dh_i.$$
Since $dy_1 \wedge dy_1 = 0$ and the cardinality of $J \cap \{1, \ldots e\}$ is at least $e - s$, we see by computing
$$\mu^*(\omega) \;=\; \mu^*(df_1) \wedge \ldots \wedge \mu^*(df_s) \wedge \mu^*(dx_J)$$
that $\mathrm{ord}_E(G) \geq (sa - 1) + (e - s)$, as stated. $\square$

## 5. Jumping numbers for graded systems

Some of the most surprising applications of multiplier ideals (e.g. [31], [7]) have involved asymptotic constructions associated to certain families of ideals. We discuss in this section the jumping coefficients that one can attach to such graded systems of ideals.

Let $X$ be a smooth complex variety. Recall from [7] that a *graded family* or *graded system* $\mathfrak{a}_\bullet = \{\mathfrak{a}_k\}_{k \in \mathbf{N}}$ of ideals is a collection of ideal sheaves $\mathfrak{a}_k \subseteq \mathcal{O}_X$ indexed by the natural numbers $\mathbf{N}$ having the property that
$$\mathfrak{a}_\ell \cdot \mathfrak{a}_m \;\subseteq\; \mathfrak{a}_{\ell+m} \quad \text{for all} \quad \ell, m \geq 1.$$
To simplify some statements, we will also assume that $\mathfrak{a}_k \neq 0$ for $k \gg 0$.

**Example 5.1** (Some graded families of ideals). (i). Fixing an ideal $\mathfrak{b} \subseteq \mathcal{O}_X$, put $\mathfrak{a}_k = \mathfrak{b}^k$. One should consider this to be a trivial example.

(ii). If $X$ is projective, let $D$ be a big divisor on $X$. Then the base-ideals
$$\mathfrak{b}_k \;=\; \mathfrak{b}\big(|kD|\big)$$
corresponding to the complete linear series $|kD|$ form a graded system. In general, one can think of graded families of ideals as essentially local objects that display the sort of complexities that arise for global linear series, and the philosophy is to use globally-inspired methods to study them.

(iii). Fixing any positive real numbers $\mu_1, \ldots, \mu_d > 0$, we get a graded family $\mathfrak{d}_\bullet$ by taking $\mathfrak{d}_k \subseteq \mathbf{C}[t_1, \ldots, t_d]$ to be the monomial ideal generated by all monomials
$$t_1^{e_1} t_2^{e_2} \cdot \ldots \cdot t_d^{e_d} \text{ with } \sum \frac{e_i}{\mu_i} \geq k.$$

When the $\mu_i$ are integral, $\mathfrak{d}_k$ is just (the integral closure of) the $k^{\text{th}}$ power of the ideal "diagonal" ideal $(t_1^{\mu_1}, \ldots, t_d^{\mu_d})$, and it is suggestive to think of $\mathfrak{d}_\bullet$ as giving meaning to the expression $(t_1^{\mu_1}, \ldots, t_d^{\mu_d})$ also when the $\mu_i$ are irrational. □

Given a graded system $\mathfrak{a}_\bullet = \{\mathfrak{a}_k\}$ and a real number $c > 0$, one can define multiplier ideals $\mathcal{J}(c \cdot \mathfrak{a}_\bullet) = \mathcal{J}(\mathfrak{a}_\bullet^c)$ by a natural asymptotic construction. These asymptotic multiplier ideals have been very important in applications.

One starts by using the Noetherian property to prove:

**Lemma 5.2.** *The ideals $\mathcal{J}(\frac{c}{p} \cdot \mathfrak{a}_p)$ all coincide for $p \gg 0$.* □

One then defines $\mathcal{J}(c \cdot \mathfrak{a}_\bullet)$ to be this common ideal, which is alternatively the unique maximal element among the ideals $\mathcal{J}(\frac{c}{p} \cdot \mathfrak{a}_p)$ appearing in the Lemma. We refer to [7], [20, Chapter 11] or [1] for details and further examples.

**Example 5.3.** For the graded family $\mathfrak{d}_\bullet$ appearing in Example 5.1(iii), $\mathcal{J}(c \cdot \mathfrak{d}_\bullet)$ is the monomial ideal spanned by all monomials $t_1^{v_1} \ldots t_d^{v_d}$ with
$$\sum \frac{v_i + 1}{\mu_i} > c. \quad □$$

The intervals of constancy for the multiplier ideals $\mathcal{J}(c \cdot \mathfrak{a}_\bullet)$ have the same shape as before:

**Lemma 5.4.** *Let $\mathfrak{a}_\bullet = \{\mathfrak{a}_k\}$ be a graded system of ideals and let $c > 0$ be a positive real number. Then*
$$\mathcal{J}(c \cdot \mathfrak{a}_\bullet) = \mathcal{J}(c' \cdot \mathfrak{a}_\bullet)$$
*for all $c' > c$ sufficiently close to $c$. In particular, the multiplier ideals $\mathcal{J}(c \cdot \mathfrak{a}_\bullet)$ are constant exactly for $c \in [\xi, \xi')$ for suitable real numbers $\xi, \xi' > 0$.*

**Definition 5.5** (Jumping coefficients for graded systems). The *jumping numbers* of $\mathfrak{a}_\bullet$ on $X$ are the endpoints of the intervals of constancy appearing in Lemma 5.4. The local jumping numbers at a fixed point $x \in X$ are defined similarly. We denote by $\text{Jump}(\mathfrak{a}_\bullet)$ the collection of all jumping coefficients of $\mathfrak{a}_\bullet$ on $X$, and by $\text{Jump}(\mathfrak{a}_\bullet; x)$ the set of all local jumping numbers. □

We shall see shortly (Example 5.10) that the collection of jumping numbers of a graded system can contain cluster points. Hence we do not try to enumerate them.

**Remark 5.6.** Lemma 5.4 may be rephrased as saying that $\text{Jump}(\mathfrak{a}_\bullet) \subseteq \mathbf{R}$ satisfies the descending chain condition (DCC): any decreasing sequence of jumping numbers stabilizes. □

**Example 5.7.** Let $\mathfrak{d}_\bullet = \{\mathfrak{d}_k\}$ be the graded family of Example 5.1(iii). Then it follows from Example 5.3 that the non-zero jumping coefficients of $\mathfrak{d}_\bullet$ consist of all real numbers of the form
$$\frac{e_1+1}{\mu_1} + \ldots + \frac{e_d+1}{\mu_d}$$
as $e_1, \ldots, e_d \in \mathbf{N}$ range over all non-negative integers. (Compare Example 1.9.) Observe that in general these are irrational and that the periodicity in Proposition 1.12 fails. $\square$

*Proof of Lemma 5.4.* Fix $p \gg 0$ which computes $\mathcal{J}(c \cdot \mathfrak{a}_\bullet)$ in the sense that $\mathcal{J}(c \cdot \mathfrak{a}_\bullet) = \mathcal{J}(\frac{c}{p} \cdot \mathfrak{a}_p)$. Applying Lemma 1.3 to $\mathfrak{a}_p$, we find that if $c' > c$ is sufficiently close to $c$, then $\mathcal{J}(\frac{c}{p} \cdot \mathfrak{a}_p) = \mathcal{J}(\frac{c'}{p} \cdot \mathfrak{a}_p)$. But $\mathcal{J}(\frac{c'}{p} \cdot \mathfrak{a}_p) \subseteq \mathcal{J}(c' \cdot \mathfrak{a}_\bullet)$, and hence
$$\mathcal{J}(c \cdot \mathfrak{a}_\bullet) \subseteq \mathcal{J}(c' \cdot \mathfrak{a}_\bullet).$$
The reverse inclusion is clear since $c' \geq c$. $\square$

Given a graded system $\mathfrak{a}_\bullet$ the log-canonical threshold $\mathrm{lct}(\mathfrak{a}_\bullet; x)$ of $\mathfrak{a}_\bullet$ at a fixed point $x \in X$ is the least real number $c > 0$ such that $\mathcal{J}(c \cdot \mathfrak{a}_\bullet)_x \neq \mathcal{O}_x X$, or equivalently the least non-zero element of $\mathrm{Jump}(\mathfrak{a}_\bullet; x)$. (If $\mathcal{J}(c \cdot \mathfrak{a}_\bullet)_x = \mathcal{O}_x X$ for all $c > 0$ we put $\mathrm{lct}(\mathfrak{a}_\bullet; x) = \infty$). The analogue of Proposition 1.17 remains valid in the current setting:

**Proposition 5.8.** *Let $\mathfrak{a}_\bullet$ be a graded system with $\mathrm{lct}(\mathfrak{a}_\bullet; x) < \infty$ for some fixed point $x$ and let $\xi \in \mathrm{Jump}(\mathfrak{a}_\bullet; x)$ be a jumping coefficient of $\mathfrak{a}_\bullet$ at $x$. Then there is another jumping number $\xi' > \xi$ of $\mathfrak{a}_\bullet$ at $x$ satisfying*
$$\xi' \leq \xi + \mathrm{lct}(\mathfrak{a}_\bullet; x).$$

*Proof.* The subadditivity theorem of [5] remains true for graded systems, i.e.
$$\mathcal{J}((c+d) \cdot \mathfrak{a}_\bullet) \subseteq \mathcal{J}(c \cdot \mathfrak{a}_\bullet) \cdot \mathcal{J}(d \cdot \mathfrak{a}_\bullet)$$
for any real numbers $c, d > 0$. This being said, the proof of Proposition 1.17 applies without change. $\square$

**Remark 5.9.** We do not know whether an analogue of Proposition 1.20 holds in the present setting. $\square$

We next give an example of a graded system whose jumping numbers contain accumulation points.

**Example 5.10** (Clustering). Working in $\mathbf{R}^2$ with coordinates $a$ and $b$, consider the region $N$ in the first quadrant given by
$$N = \left\{ (a,b) \in \mathbf{R}^2_{\geq 0} \,\big|\, (a-1)(b-1) \geq 1 \right\}$$
(Figure 2). Inspired by a construction of Mustaţă [29], the idea is to realize $N$ as the the limit of the Newton polyhedra coming from a graded sequence of monomial ideals. Choose a nested sequence of closed rational convex polyhedra
$$N_1 \subseteq N_2 \subseteq N_3 \subseteq \ldots$$

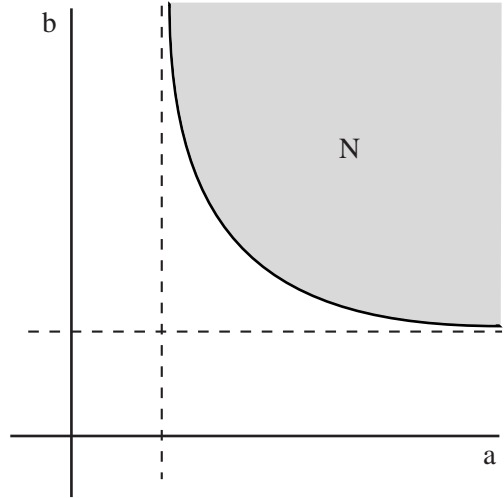

FIGURE 2. Graded system of ideals with clustering jumping numbers

contained in $N$, with $N_k + \mathbf{R}_{\geq 0}^2 \subseteq N_k$ for each $k$, approximating $N$ in the sense that $N = \cup N_k$. Put $P_k = k \cdot N_k$. Then the $P_k$ form a graded family with respect to Minkowski addition, i.e.
$$P_\ell + P_m \subseteq P_{\ell+m} \quad \text{for all} \quad \ell, m \geq 1.$$
Let $\mathfrak{a}_k \subseteq \mathbf{C}[s,t]$ be the monomial ideal spanned by all monomials whose exponent vectors lie in $P_k$. Then $\mathfrak{a}_\bullet = \{\mathfrak{a}_k\}$ is a graded system of monomial ideals.

The multiplier ideals of $\mathfrak{a}_\bullet$ are computed from $N$ in much the same way that the multiplier ideals of a single monomial ideal are computed from its Newton polyhedron. Specifically, it follows from Example 1.8 and the constructions that $\mathcal{J}(c \cdot \mathfrak{a}_\bullet)$ is the monomial ideal spanned by all monomials $s^e t^f$ such that
$$(e+1, f+1) \in \text{interior}(c \cdot N).$$
Therefore each lattice point $(e, f) \in \mathbf{N}^2$ in the first quadrant determines a jumping number $\xi = \xi_{(e,f)}$ by the condition that $(e+1, f+1)$ lies on the boundary of $\xi \cdot N$, or in other words that
$$\left(\frac{e+1}{\xi} - 1\right)\left(\frac{f+1}{\xi} - 1\right) = 1.$$
Thus the jumping numbers of $\mathfrak{a}_\bullet$ consist precisely of the rational numbers
$$\left\{\frac{(e+1)(f+1)}{e+f+2} \,\Big|\, e, f \geq 0\right\},$$
and by fixing $e$ and letting $f \to \infty$ one sees that these cluster at all positive integers.[8] □

**Remark 5.11.** Note that each of the ideals $\mathfrak{a}_k$ constructed in the previous example vanishes along the two coordinate axes in $\mathbf{C}^2$. We will see momentarily that if $\mathfrak{a}_\bullet = \{\mathfrak{a}_k\}$ is a graded family of ideals which vanish only on a finite set, then $\text{Jump}(\mathfrak{a}_\bullet)$ has no accumulation points. $\square$

---

[8]It is amusing to verify explicitly that all positive integers actually occur among the stated jumping numbers, as required by Lemma 5.4!

Suppose now that $\mathfrak{a}_\bullet = \{\mathfrak{a}_k\}$ is a graded system of ideals having the property that each $\mathfrak{a}_k$ vanishes only at a single point $x \in X$. Denote by $\mathfrak{m} \subseteq \mathcal{O}_X$ the maximal ideal of $x$.

**Lemma 5.12.** *The collection $\text{Jump}(\mathfrak{a}_\bullet) = \text{Jump}(\mathfrak{a}_\bullet; x)$ of jumping numbers of the $\mathfrak{m}$-primary family $\mathfrak{a}_\bullet$ is discrete, and each of the multiplier ideals $\mathcal{J}(c \cdot \mathfrak{a}_\bullet)$ is $\mathfrak{m}$-primary, hence of finite codimension in $\mathcal{O}_x X$.*

Just as in Definition 1.22 one can then attach a multiplicity to each jumping coefficient, and we write

$$(31) \qquad 0 = \kappa_0(\mathfrak{a}_\bullet) \leq \kappa_1(\mathfrak{a}_\bullet) \leq \kappa_2(\mathfrak{a}_\bullet) \leq \ldots$$

for the sequence of jumping numbers of $\mathfrak{a}_\bullet$ at $x$, each repeated according to its multiplicity.

*Proof of Lemma 5.12.* For any integer $\ell > 0$ there is an inclusion $\mathfrak{a}_\ell \subseteq \mathcal{J}(\ell \cdot \mathfrak{a}_\bullet)$. Therefore $\mathcal{J}(\ell \cdot \mathfrak{a}_\bullet)$ is $\mathfrak{m}$-primary, and the number of jumping numbers of $\mathfrak{a}_\bullet$ which are $\leq \ell$ is bounded by the codimension of $\mathfrak{a}_\ell$. In particular, $\text{Jump}(\mathfrak{a}_\bullet)$ is discrete. $\square$

Our next object is to study the variational properties of these jumping numbers. We start with a general lemma about restrictions of asymptotic multiplier ideals.

**Lemma 5.13.** *Let $Y$ be a smooth variety, and $f : Y \longrightarrow T$ a smooth mapping from $Y$ onto a smooth variety $T$. Fix a graded system $\mathfrak{a}_\bullet = \{\mathfrak{a}_k\}$ of ideals on $Y$. Given $t \in T$ write $Y_t$ for the fibre of $Y$ over $t \in T$, and denote by $\mathfrak{a}_{\bullet,t}$ the graded family of ideals obtained by restricting $\mathfrak{a}_\bullet$ to $Y_t$. Assume that none of the $\mathfrak{a}_k$ vanish on any of the fibres $Y_t$.*

(i). *For arbitrary $t \in T$ and $c > 0$ there is an inclusion*

$$(32) \qquad \mathcal{J}(Y_t, c \cdot \mathfrak{a}_{\bullet,t}) \subseteq \mathcal{J}(Y, c \cdot \mathfrak{a}_\bullet) \cdot \mathcal{O}_{Y_t}.$$

(ii). *There is a countable union of proper closed subvarieties $\mathcal{B} \subsetneq T$ such that if $t \in T - \mathcal{B}$ then equality holds in (32) for every $c > 0$.*

*Proof.* Given $t \in T$ and $c > 0$, fix $p \gg 0$ so that

$$\mathcal{J}(Y, c \cdot \mathfrak{a}_\bullet) = \mathcal{J}(Y, \tfrac{c}{p} \cdot \mathfrak{a}_p) \text{ and } \mathcal{J}(Y_t, c \cdot \mathfrak{a}_{\bullet,t}) = \mathcal{J}(Y_t, \tfrac{c}{p} \cdot \mathfrak{a}_{p,t}).$$

The restriction theorem implies that $\mathcal{J}(Y_t, \tfrac{c}{p} \cdot \mathfrak{a}_{p,t}) \subseteq \mathcal{J}(Y, \tfrac{c}{p} \cdot \mathfrak{a}_p) \cdot \mathcal{O}_{Y_t}$, and (i) follows. For (ii), fix any $p > 0$. By the theorem on generic restrictions of multiplier ideals, there exists a proper closed subset $B_p \subsetneq T$ such that

$$\mathcal{J}(Y, d \cdot \mathfrak{a}_p) \cdot \mathcal{O}_{Y_t} = \mathcal{J}(Y_t, d \cdot \mathfrak{a}_{p,t}) \quad \text{for all } t \in T - B_p \text{ and all } d > 0.$$

Take $\mathcal{B} = \cup_{p \geq 1} B_p$, and fix $c > 0$. There is a natural number $p \gg 0$ (depending on $c$) such that $\mathcal{J}(Y, c \cdot \mathfrak{a}_\bullet) = \mathcal{J}(Y, \tfrac{c}{p} \cdot \mathfrak{a}_p)$. If $t \in T - \mathcal{B}$ then by construction $\mathcal{J}(Y, \tfrac{c}{p} \cdot \mathfrak{a}_p)|Y_t =$

$\mathcal{J}(Y_t, \frac{c}{p} \cdot \mathfrak{a}_{p,t})$, and so

$$\mathcal{J}(Y, c \cdot \mathfrak{a}_\bullet) \cdot \mathcal{O}_{Y_t} \;=\; \mathcal{J}(Y_t, \tfrac{c}{p} \cdot \mathfrak{a}_{p,t}) \;\subseteq\; \mathcal{J}(Y_t, c \cdot \mathfrak{a}_{\bullet,t}).$$

On the other hand, we have the reverse inclusion from (i) and consequently $\mathcal{J}(Y, c \cdot \mathfrak{a}_\bullet)|Y_t = \mathcal{J}(Y_t, c \cdot \mathfrak{a}_{\bullet,t})$, as required. □

One then has an analogue of Proposition 1.24.

**Proposition 5.14** (Semicontinuity for graded systems). *Let $T$ be a smooth curve, and let $\mathfrak{a}_\bullet = \{\mathfrak{a}_k\}$ be graded family of ideals on $X \times T$, each of whose zeroes are supported on $\{x\} \times T$ for some $x \in X$. Given $t \in T$ write $\mathfrak{a}_{\bullet,t}$ for the graded system on $X \times \{t\} = X$ arising from the restriction of $\mathfrak{a}_\bullet$. Then there is a countable set $\mathcal{B} \subseteq T$ such that all the jumping numbers $\kappa_i(\mathfrak{a}_{\bullet,t})$ are constant for $t \in T - \mathcal{B}$. Moreover if $t^* \in T$ is an arbitrary point and $t \in T - \mathcal{B}$, then*

(33) $$\kappa_i(\mathfrak{a}_{\bullet,t^*}) \;\leq\; \kappa_i(\mathfrak{a}_{\bullet,t}) \quad \text{for all} \quad i \geq 0.$$

*Remarks on the proof.* The proof of 1.24 goes through with only minor changes using Lemma 5.13 in place of the restriction and generic restriction theorems. We leave details to the reader. □

Finally, we make a few remarks about jumping coefficients associated to plurisubharmonic functions. Let $X$ be a complex manifold, and let $\varphi$ be a plurisubharmonic (or PSH) function on $X$. The multiplier ideal sheaf $\mathcal{J}(X, \varphi)$ is defined to be the analytic sheaf of ideals whose stalk at $x \in X$ consists of all germs of holomorphic functions $f$ such that $|f|^2 e^{-\varphi}$ is integrable in some neighborhood of $x$ in $X$. If $f_1, ..., f_N$ are analytic functions on $X$, then $\varphi = c \cdot \log(|f_1|^2 + ... + |f_N|^2)$ is psh, and thus the notion of multiplier ideal associated to a PSH function generalizes the multiplier ideals assoicated to effective **Q**-divisors and to ideal sheaves as in Remark 1.2.

Given a PSH function $\varphi$ and a fixed point $x \in X$, there is a natural way to define the jumping numbers of $\varphi$ at $x$. Specifically, a positive real number $\xi > 0$ is a jumping coefficient of $\varphi$ at $x$ if

$$\mathcal{J}(c \cdot \varphi)_x \;\subsetneq\; \mathcal{J}(c' \cdot \varphi)_x$$

whenever $c > \xi > c'$. We denote by $\text{Jump}(\varphi; x)$ the collection of all jumping numbers at $x$. Some of the properties discussed above — for example Proposition 5.8 — extend without change. However others seems more subtle. In particular, it appears difficult to determine whether or not the analogue of Lemma 5.4 remains valid for PSH functions, i.e. whether $\text{Jump}(\varphi; x)$ satisfies the DCC. When $\dim X = 2$, Favre and Jonsson [11] have made some important progress on this question, but their methods do not seem to generalize to higher dimensions.

Department of Mathematics, University of Illinois at Chicago
851 South Morgan Street (M/C 249), Chicago, IL 60607-7045, USA

*E-mail address*: `ein@math.uic.edu`

Department of Mathematics, University of Michigan, Ann Arbor, MI 48109, USA

*E-mail address*: `rlaz@umich.edu`

Department of Mathematics, University of Michigan, Ann Arbor, MI 48109

*E-mail address*: `kesmith@umich.edu`

Department of Mathematics, University of Illinois at Urbana-Champaign, Urbana, IL 61801

*E-mail address*: `dror@math.uiuc.edu`